\documentclass{amsart}

\usepackage{amssymb}
\usepackage{fancyhdr}
\usepackage{amsmath}
\usepackage{mathrsfs}
\usepackage{amsthm}
\usepackage{bm}
\usepackage{verbatim}
\usepackage{booktabs}
\usepackage{dcolumn}
\usepackage{pgf}
\usepackage{float}
\usepackage[square,numbers]{natbib}
\usepackage{xxcolor}
\usepackage[colorlinks=true, linkcolor=blue, urlcolor=blue, citecolor = blue]{hyperref}
\usepackage{stmaryrd}
\usepackage[squaren]{SIunits}
\usepackage{graphicx}
\usepackage{tabularx}
\usepackage{chngcntr}
\usepackage{subcaption}
\usepackage{afterpage}
\usepackage{multirow}
\usepackage[normalem]{ulem}

\newtheorem{theorem}{Theorem}[section]

\theoremstyle{definition}

\theoremstyle{remark}

\numberwithin{equation}{section}

\def\dd{\, \mathrm{d}}

\newtheorem{innercustomgeneric}{\customgenericname}
\providecommand{\customgenericname}{}
\newcommand{\newcustomproblem}[2]{%
	\newenvironment{#1}[1]
	{%
		\renewcommand\customgenericname{#2}%
		\renewcommand\theinnercustomgeneric{##1}%
		\innercustomgeneric
	}
	{\endinnercustomgeneric}
}

\newcustomproblem{customprob}{Problem}
\newcustomproblem{customalgo}{Algorithm}
\newcommand*{\bqed}{\hfill\ensuremath{\blacksquare}}%

\begin{document}
	
	
	\title[Numerical validation of Koiter's model]{Numerical validation of Koiter's model for all the main types of linearly elastic shells in the static case}
	
	
	\author{Wangxi Duan}
	\address{Department of Applied Mathematics, School of Sciences, Xi'an University of Technology, P.O.Box 1243,Yanxiang Road NO.58 XI'AN, 710054, Shaanxi Province, China}
	\email{duan\_wx97@163.com}
	
	\author{Paolo Piersanti*}
	\address{Department of Mathematics and Institute for Scientific Computing and Applied Mathematics, Indiana University Bloomington, 729 East Third Street, Bloomington, Indiana, USA}
	\email{ppiersan@iu.edu}
	
	\author{Xiaoqin Shen*}
	\address{Department of Applied Mathematics, School of Sciences, Xi'an University of Technology, P.O.Box 1243,Yanxiang Road NO.58 XI'AN, 710054, Shaanxi Province, China}
	\email{xqshen@xaut.edu.cn}
	
	\author{Qian Yang}
	\address{Department of Applied Mathematics, School of Sciences, Xi'an University of Technology, P.O.Box 1243,Yanxiang Road NO.58 XI'AN, 710054, Shaanxi Province, China}
	\email{yq931122@sina.com}
	
\begin{abstract}
In this paper, we validate the main theorems establishing the justification of Koiter's model, established by Ciarlet and his associates, for all the types of linearly elastic shells via a set of numerical experiments.
\end{abstract}

\maketitle
\renewcommand{\thefootnote}{\fnsymbol{footnote}}
\footnotetext[1]{Corresponding author.}

\tableofcontents

\section{Introduction} \label{Sec:0}

The modelling of a thin linearly elastic shell can be \emph{either} performed by resorting to the classical equations of the three-dimensional linearised elasticity, \emph{or} by means of a specific set two-dimensional equations obtained by W.T. Koiter in the seminal works~\cite{Koiter,Koiter1970}.

By means of a rigorous asymptotic analysis (cf.~\cite{Ciarlet2000} and~\cite{Ciarlet19962}) it is possible to prove that both the three-dimensional equations of linearised elasticity and the two-dimensional equations of Koiter's model have the same asymptotic behaviour as the thickness of the shell tend to zero. Moreover, the striking feature of Koiter's model is that it does not \emph{a priori} depend on the geometrical features of the shell under consideration, as it indeed includes both the \emph{elliptic membrane shells} bilinear form and the \emph{flexural shells} bilinear form.

This feature is very convenient if one aims to study a variety of situations that would be too complicated to address by means of the classical three-dimensional approach, like, for instance when the shell is anisotropic, inhomogeneous (cf., e.g., \cite{CaillerieSanchez1995a} and~\cite{CaillerieSanchez1995b}), or or its thickness varies periodically (cf., e.g., \cite{TelegaLewinski1998a} and~\cite{TelegaLewinski1998b}).

Koiter's model has been widely studied by Ciarlet and his group~\cite{CiaLods1996a,CiaLods1996b,Ciarlet19962,CiaPie2018b,CiaPie2018bCR, Ciarlet19963,CiaMarPie2018b,CiaMarPie2018}. For what concerns the literature related to time-dependent problems in elasticity, it is worth mentioning the papers~\cite{BockJar2009,BockJar2013,BockJarSil2016,Pie2020,Pie2020-1, Pie2021,Shen2018}. It is also worth mentioning the papers~\cite{Pie2020-2,PS,Rodri2018,Shen2019,ShenPiePie2020,Shen2021,Shen2020}.

The purpose of this paper is to validate the main theorems establishing the justification of Koiter's model, established by Ciarlet and his associates, for all the types of linearly elastic shells via a set of numerical experiments. More precisely, we present a set of numerical experiments showing that, for all the main types of linearly elastic shells (i.e., linearly elastic elliptic membrane shells, linearly elastic generalized membrane shells and linearly elastic flexural shells) the solution of Koiter's model asymptotically behaves as the solution of the model based on the three-dimensional equations of linearized elasticity as the thickness approaches zero.

The motivations that led us to carry out the present investigation are that, despite Koiter's model has been rigorously proved from the mathematical point of view, no numerical simulations have been provided in this direction. The availability of numerical results would allow us to gain more insight into the \emph{quantitative} way the solution of Koiter's model asymptotically behaves like the original three-dimensional model.

This paper is divided into seven sections (including this one). In section~\ref{Sec:1} we recall the notation and some geometrical preliminaries we will be making use of throughout the paper; in section~\ref{Sec:2} we recall the formulation of the three-dimensional problem based on the classical equations of linearized elasticity; in section~\ref{Koiter-model} we recall the definition and main properties of Koiter's model; in section~\ref{Sec:3} we recall the main result ensuring the justification of Koiter's model for linearly elastic elliptic membrane shells and we present some numerical experiments validating the converegences established in the aforementioned theorem.
In section~\ref{Sec:4} and~\ref{Sec:5} we repeat the same type of investigation for generalized membrane shells of the first kind and for flexural shells respectively. Finally, in section~\ref{Sec:6} we draw a summary of the numerical results we presented and we talk about future projects.

\section{Geometrical preliminaries}
\label{Sec:1}

The geometrical notation is borrowed from the reference~\cite{Ciarlet2000} (see also \cite{Ciarlet2005}).

Greek indices, except the thickness parameter $\varepsilon$, take their values in the set $\{1,2\}$. Latin indices, instead, take their values in the set $\{1,2,3\}$. The Einstein summation convention with respect to repeated indices is systematically used in conjunction with these two rules. The symbol $\mathbb{E}^3$ indicates the three-dimensional Euclidean space whose origin is denoted by $\bm{O}$; the Euclidean inner product and the vector product of $\bm{u}, \bm{v} \in \mathbb{E}^3$ are denoted by $\bm{u} \cdot \bm{v}$ and $\bm{u} \times \bm{v}$; the Euclidean norm of $\bm{u} \in \mathbb{E}^3$ is denoted by $\left|\bm{u} \right|$. The notation $\delta^j_i$ designates the standard Kronecker symbol.

Given an open subset $\Omega$ of $\mathbb{R}^n$, notations such as $L^2(\Omega)$, $H^1(\Omega)$, or $H^2 (\Omega)$, designate the usual Lebesgue and Sobolev spaces. The notation $\left\| \cdot \right\|_X$ designates the norm in a normed vector space $X$. Spaces of vector-valued functions are denoted with boldface letters.

A \emph{domain in} $\mathbb{R}^n$ is a bounded and connected open subset $\Omega$ of $\mathbb{R}^n$, whose boundary $\partial \Omega$ is Lipschitz-continuous, the set $\Omega$ being locally on a single side of $\partial \Omega$.

Let $\omega$ be a domain in $\mathbb{R}^2$, and let us indicate a generic point in $\omega$ by $y = (y_\alpha)$. A mapping $\bm{\theta} \in \mathcal{C}^1(\overline{\omega}; \mathbb{E}^3)$ is said to be an \emph{immersion} if the two vectors
\[
\bm{a}_\alpha (y) := \partial_\alpha \bm{\theta} (y)
\]
are linearly independent at each point $y \in \overline{\omega}$. Therefore, the set $\bm{\theta}(\overline{\omega})$, namely the image of the set $\overline{\omega}$ under the mapping $\bm{\theta}$, is a \emph{surface in} $\mathbb{E}^3$ equipped with $y_1, y_2$ as its \emph{curvilinear coordinates}. 
For each $y\in \overline{\omega}$, the vectors $\bm{a}_\alpha (y)$ span the plane which is tangential to the surface $\bm{\theta} (\overline{\omega})$ at the point $\bm{\theta} (y)$.

Define the unit vector
\[
\bm{a}_3 (y) := \frac{\bm{a}_1(y) \wedge \bm{a}_2 (y)}{|\bm{a}_1(y) \wedge \bm{a}_2 (y)|},
\]
that it is normal to $\bm{\theta}(\overline{\omega})$ at the point $\bm{\theta}(y)$. Define the vectors $\bm{a}^j(y)$ by
\[
\bm{a}^j(y) \cdot \bm{a}_i (y) = \delta^j_i.
\]

The vectors $\bm{a}^i$ form the \emph{contravariant} basis at $\bm{\theta} (y)$.

The first fundamental form of the surface $\bm{\theta} (\overline{\omega})$ can be defined by means of its covariant components in the following fashion:
\[
a_{\alpha \beta} := \bm{a}_\alpha \cdot \bm{a}_\beta = a_{\beta \alpha} \in \mathcal{C}^0 (\overline{\omega}).
\]

The \emph{contravariant components} of the first fundamental form are, instead, defined as follows:
\[
a^{\alpha \beta}:= \bm{a}^\alpha \cdot \bm{a}^\beta = a^{\beta \alpha}\in \mathcal{C}^0(\overline{\omega}).
\]

The contravariant components of the first fundamental form of the surface $\bm{\theta}(\overline{\omega})$ all together define a matrix field $(a^{\alpha \beta})$, which can be proved to be symmetric and equal to the inverse of the matrix field $(a_{\alpha \beta})$, itself defined in terms of the covariant components of the first fundamental form of the surface $\bm{\theta}(\overline{\omega})$. 

Moreover, we have that $\bm{a}^\beta = a^{\alpha \beta}\bm{a}_\alpha$ and $\bm{a}_\alpha = a_{\alpha \beta} \bm{a}^\beta$. The \emph{surface element} along $\bm{\theta} (\overline{\omega})$ is defined at each point $\bm{\theta} (y), \, y \in \overline{\omega}$, by $\sqrt{a(y)}\dd y$, where
\[
a := \det (a_{\alpha \beta}) \in \mathcal{C}^0 (\overline{\omega}).
\]

Let $\bm{\theta} \in \mathcal{C}^2(\overline{\omega}; \mathbb{E}^3)$ be an immersion. We can define the \emph{second fundamental form} of the surface $\bm{\theta}(\overline{\omega})$ in terms of its \emph{covariant components} as follows:
\[
b_{\alpha \beta}:= \partial_\alpha \bm{a}_\beta \cdot \bm{a}_3 = - \bm{a}_\beta \cdot \partial_\alpha \bm{a}_3 = b_{\beta \alpha} \in \mathcal{C}^0(\overline{\omega}).
\]

Likewise, it is possible define the second fundamental form of the surface $\bm{\theta}(\overline{\omega})$ in terms of its \emph{mixed components} as follows:
\[
b^\beta_\alpha := a^{\beta \sigma} b_{\alpha \sigma} \in \mathcal{C}^0(\overline{\omega}).
\]

Define the \emph{Christoffel symbols} corresponding to the immersion $\bm{\theta}$ by the following formulas:
\[
\Gamma^\sigma_{\alpha \beta}:= \partial_\alpha \bm{a}_\beta \cdot \bm{a}^\sigma = \Gamma^\sigma_{\beta \alpha} \in \mathcal{C}^0 (\overline{\omega}).
\]

The classification of the shell models we will be considering in the forthcoming sections somehow hinges on the \emph{Gaussian curvature} of the surface $\bm{\theta} (\overline{\omega})$ at each point $\bm{\theta}(y) , \, y \in \overline{\omega}$, which is defined by:
\[
\kappa (y) := \frac{\det (b_{\alpha \beta} (y))}{\det (a_{\alpha \beta} (y))} = \det \left( b^\beta_\alpha (y)\right).
\]

We can associate with the displacement field $\eta_i \bm{a}^i$ the \emph{linearized change of metric tensor}, expressed in terms of its covariant components:
\[
\gamma_{\alpha \beta}(\bm{\eta}) := \frac12 (\partial_\beta \eta_\alpha + \partial_\alpha \eta_\beta ) - \Gamma^\sigma_{\alpha \beta} \eta_\sigma - b_{\alpha \beta} \eta_3 = \gamma_{\beta \alpha} (\bm{\eta}).
\]

Likewise, we can associate with the displacement field $\eta_i \bm{a}^i$ the \emph{linearized change of curvature tensor}, expressed in terms of its covariant components:
\begin{align*}
	\rho_{\alpha \beta}(\bm{\eta}) &:= \partial_{\alpha \beta}\eta_3 -\Gamma_{\alpha\beta}^\sigma \partial_\sigma\eta_3 -b_\alpha^\sigma b_{\sigma\beta}\eta_3\\
	&+b_\alpha^\sigma(\partial_\beta \eta_\sigma-\Gamma_{\beta\sigma}^\tau\eta_\tau)+b_\beta^\tau(\partial_\alpha\eta_\tau-\Gamma_{\alpha\tau}^\sigma\eta_\sigma)\\
	&+(\partial_\alpha b_\beta^\tau+\Gamma_{\alpha\sigma}^\tau b_\beta^\sigma-\Gamma_{\alpha\beta}^\sigma b_\sigma^\tau)\eta_\tau=\rho_{\beta\alpha}(\bm{\eta}).
\end{align*}

\section{The three-dimensional problem for a ``general'' linearly elastic shell} \label{Sec:2}

Let $\omega$ be a domain in $\mathbb{R}^2$, let $\gamma:= \partial \omega$, and let $\gamma_0$ be a non-empty relatively open subset of $\gamma$. For each $\varepsilon > 0$, define the sets
$$
\Omega^\varepsilon := \omega \times \left] - \varepsilon , \varepsilon \right[ \textup{ and } \Gamma^\varepsilon_0 := \gamma_0 \times \left[ - \varepsilon , \varepsilon \right] \textup{ and } \Gamma^\varepsilon_\pm := \omega \times \{\pm \varepsilon\},
$$
we let $x^\varepsilon = (x^\varepsilon_i)$ denote a generic point in the set $\overline{\Omega^\varepsilon}$, and we define $\partial^\varepsilon_i := \partial / \partial x^\varepsilon_i$. As a result, we also have $x^\varepsilon_\alpha = y_\alpha$ and $\partial^\varepsilon_\alpha = \partial_\alpha$.

Let $\bm{\theta} \in \mathcal{C}^3(\overline{\omega}; \mathbb{E}^3)$ be an immersion and let $\varepsilon > 0$. Consider a linearly elastic shell with \emph{middle surface} $\bm{\theta} (\overline{\omega})$ and with \emph{constant thickness} $2 \varepsilon$. The \emph{reference configuration} of one such linearly elastic shell is thus identified by the set $\bm{\Theta} (\overline{\Omega^\varepsilon})$, where the mapping $\bm{\Theta} : \overline{\Omega^\varepsilon} \to \mathbb{E}^3$ is defined by:
\[
\bm{\Theta} (x^\varepsilon) := \bm{\theta} (y) + x^\varepsilon_3 \bm{a}^3(y) \text{ at each point } x^\varepsilon = (y, x^\varepsilon_3) \in \overline{\Omega^\varepsilon}.
\]

Associated with the immersion $\bm{\Theta}$ is the following metric tensor, defined by means of its \emph{covariant components}:
\[
g^\varepsilon_{ij} := \bm{g}^\varepsilon_i \cdot \bm{g}^\varepsilon_j \in \mathcal{C}^1(\overline{\Omega^\varepsilon}).
\]

The same metric tensor can be defined by means of its \emph{contravariant components} in the following fashion:
\[
g^{ij, \varepsilon} := \bm{g}^{i, \varepsilon} \cdot \bm{g}^{i,\varepsilon} \in \mathcal{C}^1(\overline{\Omega^\varepsilon}).
\]

The contravariant components of the aforementioned metric tensor define matrix field $(g^{ij, \varepsilon})$ which can be proved to be symmetric and equal to the inverse of the matrix field $(g^\varepsilon_{ij})$.

It can also be proved that the following identities hold:
\begin{align*}
\bm{g}^{j, \varepsilon} &= g^{ij, \varepsilon} \bm{g}^\varepsilon_i,\\
g^\varepsilon_i &= g^\varepsilon_{ij} \bm{g}^{j, \varepsilon}.
\end{align*}

The \emph{volume element} in $\bm{\Theta} (\overline{\Omega^\varepsilon})$ at each point $\bm{\Theta} (x^\varepsilon)$, $x^\varepsilon \in \overline{\Omega^\varepsilon}$, is defined by $\sqrt{g^\varepsilon (x^\varepsilon)} \dd x^\varepsilon$, where:
\[
g^\varepsilon := \det (g^\varepsilon_{ij}) \in \mathcal{C}^1(\overline{\Omega^\varepsilon}).
\]

Associated with the immersion $\bm{\Theta}$ are the \emph{Christoffel symbols} defined by:
\[
\Gamma^{p, \varepsilon}_{ij}:= \partial_i \bm{g}^\varepsilon_j \cdot \bm{g}^{p, \varepsilon} = \Gamma^{p, \varepsilon}_{ji} \in \mathcal{C}^0(\overline{\Omega^\varepsilon}).
\]

It is worth recalling (cf., e.g., \cite{Ciarlet2000}) that $\Gamma^{3,\varepsilon}_{\alpha 3} = \Gamma^{p, \varepsilon}_{33} = 0$.

Associated with the displacement field $v^\varepsilon_i \bm{g}^{i, \varepsilon}$ is the \emph{linearized strain tensor} defined by means of its covariant components in the following fashion:
\[
e^\varepsilon_{i\|j} (\bm{v}^\varepsilon) :=\frac12 (\partial^\varepsilon_j v^\varepsilon_i + \partial^\varepsilon_i v^\varepsilon_j) - \Gamma^{p, \varepsilon}_{ij} v^\varepsilon_p = e_{j\|i}^\varepsilon (\bm{v}^\varepsilon).
\]

The functions $e^\varepsilon_{i\|j} (\bm{v}^\varepsilon)$ are customarily referred to as the \emph{linearized strains in curvilinear coordinates} associated with the displacement field $v^\varepsilon_i \bm{g}^{i, \varepsilon}$.

Throughout the rest of this paper it is assumed that the material constituting the shell is \emph{homogeneous}, \emph{isotropic}, and \emph{linearly elastic}. The rheology of an elastic material with these mechanical features is entirely governed by its two \emph{Lam\'{e} constants} $\lambda \ge 0$ and $\mu > 0$ (for details, see, e.g., Section 3.8 of~\cite{Ciarlet1988}).

The linearly elastic shells we shall be considering are all subjected to \emph{applied body forces}, whose density per unit volume is defined by means of their covariant components $f^{i, \varepsilon} \in L^2(\Omega^\varepsilon)$, and to \emph{applied surface forces} whose density per unit area is defined by means of their covariant components $h^{i, \varepsilon} \in L^2(\Gamma^\varepsilon_{+} \cup \Gamma^\varepsilon_{-})$.
Finally, for sake of well-posedness purposes, we need to assume that the shell under consideration is also subjected to a \emph{homogeneous boundary condition of place} along the portion $\Gamma^\varepsilon_0$ of its lateral face (i.e., the displacement vanishes on $\Gamma^\varepsilon_0$).

The problem of finding the equilibrium position of a three-dimensional linearly elastic shell on which a certain deforming load acts is classical, and amounts to minimizing an \emph{ad hoc} energy functional, that will be denoted by $J^\varepsilon$ in what follows.

The functions
\[
A^{ijk\ell, \varepsilon} := \lambda g^{ij, \varepsilon} g^{k\ell, \varepsilon} + \mu \left( g^{ik, \varepsilon} g^{j\ell, \varepsilon} + g^{i\ell, \varepsilon} g^{jk, \varepsilon} \right) =
A^{jik\ell, \varepsilon} =  A^{k\ell ij, \varepsilon},
\]
denote the contravariant components of the \emph{elasticity tensor} of the linearly elastic material constituting the shell. Then the unknown of the problem, which is the displacement vector field $\bm{u}^\varepsilon = (u^\varepsilon_i)$. The to-be-minimized \emph{energy} functional $J^\varepsilon : \bm{H}^1(\Omega^\varepsilon) \to \mathbb{R}$ is defined by:
\[
J^\varepsilon (\bm{v}^\varepsilon) := \frac12 \int_{\Omega^\varepsilon} A^{ijk\ell, \varepsilon} e^\varepsilon_{k\| \ell}  (\bm{v}^\varepsilon)e^\varepsilon_{i\|j} (\bm{v}^\varepsilon) \sqrt{g^\varepsilon} \dd x^\varepsilon - \int_{\Omega^\varepsilon} f^{i, \varepsilon} v^\varepsilon_i \sqrt{g^\varepsilon} \dd x^\varepsilon,
\]
for each $\bm{v}^\varepsilon = (v^\varepsilon_i) \in \bm{V}(\Omega^\varepsilon) := \{\bm{v}^\varepsilon = (v^\varepsilon_i) \in \bm{H}^1(\Omega^\varepsilon) ; \; \bm{v}^\varepsilon = \textbf{0} \text{ on } \Gamma^\varepsilon_0\}$.

This \emph{minimization problem} admits a unique solutions, and is equivalent to the following set of variational inequalities:
\begin{customprob}{$\mathcal{P}(\Omega^\varepsilon)$}\label{problem0}
	Find $\bm{u}^\varepsilon \in \bm{V}(\Omega^\varepsilon)$ that satisfies the following variational equations:
	$$
	\int_{\Omega^\varepsilon}
	A^{ijk\ell, \varepsilon} e^\varepsilon_{k\| \ell} (\bm{u}^\varepsilon) e^\varepsilon_{i\| j}(\bm{v}^\varepsilon) \sqrt{g^\varepsilon} \dd x^\varepsilon = \int_{\Omega^\varepsilon} f^{i , \varepsilon} v^\varepsilon_i\sqrt{g^\varepsilon} \dd x^\varepsilon + \int_{\Gamma^\varepsilon_{+} \cup \Gamma^\varepsilon_{-}} h^{i,\varepsilon} v_i \sqrt{g(\varepsilon)} \dd \Gamma^\varepsilon,
	$$
	for all $\bm{v}^\varepsilon = (v^\varepsilon_i) \in \bm{V}(\Omega^\varepsilon)$.
	\bqed	
\end{customprob}

It was shown in the seminal papers~\cite{Ciarlet1996,Ciarlet19962,Ciarlet19964,Ciarlet19963} that if the shell thickness $\varepsilon$ is let approach zero in the variational equations of Problem~\ref{problem0}, we obtain a set of two-dimensional equations (\emph{two-dimensional} in the sense that they are defined over a domain in $\mathbb{R}^2$). The form of this set of two-dimensional variational equations solely depends on the assumptions on the data we made beforehand and the geometrical properties of the middle surface of the shell under considerations. 

Our objective is to verify, via \emph{ad hoc} numerical simulations, that under the assumptions on the data set out in~\cite{Ciarlet1996,Ciarlet19962,Ciarlet19964,Ciarlet19963} the two-dimensional models recovered in these papers are indeed the correct ones.

\section{Koiter's model}
\label{Koiter-model}

The model announced in Progblem~\ref{problem0}, that is governed by the three-dimensional equations of linearised elasticity, comes with a series of drawbacks that cannot in general be treated via the standard techniques. For this reason, the set of \emph{two-dimensional}equations devised in 1970 by W.T. Koiter in the seminal paper~\cite{Koiter1970} are preferable for modelling the displacement of a linearly elastic shell (``two-dimensional'' in the sense that it is posed over $\overline{\omega}$ instead of $\overline{\Omega^\varepsilon}$) that could not otherwise be treated via the standard approach.

The functional space where Koiter's model is posed is the following
$$
\bm{V}_K (\omega):= \{\bm{\eta}=(\eta_i) \in H^1(\omega)\times H^1(\omega)\times H^2(\omega);\eta_i=\partial_{\nu}\eta_3=0 \textup{ on }\gamma_0\},
$$
where the symbol $\partial_{\nu}$ denotes the \emph{outer unit normal derivative operator along} $\gamma$, and the subscript $K$ aptly recalls the connection with Koiter's model. Equip the space $\bm{V}_K(\omega)$ the norm $\|\cdot\|_{\bm{V}_K (\omega)}$ by
$$
\|\bm{\eta}\|_{\bm{V}_K(\omega)}:=\left\{\sum_{\alpha}\|\eta_\alpha\|_{1,\omega}^2+\|\eta_3\|_{2,\omega}^2\right\}^{1/2}\quad\mbox{ for each }\bm{\eta}=(\eta_i)\in \bm{V}_K(\omega).
$$

Recall that the \emph{fourth-order two-dimensional elasticity tensor of the shell}, viewed here as a two-dimensional linearly elastic body, is expressed by in terms of its contravariant components in the following fashion:
$$
a^{\alpha \beta \sigma \tau} = \frac{4\lambda \mu}{\lambda + 2 \mu} a^{\alpha \beta} a^{\sigma \tau} + 2\mu \left(a^{\alpha \sigma} a^{\beta \tau} + a^{\alpha \tau} a^{\beta \sigma}\right).
$$

Finally, define the bilinear forms $B_M(\cdot,\cdot)$ and $B_F(\cdot,\cdot)$ by
\begin{align*}
	B_M(\bm{\xi}, \bm{\eta})&:=\int_\omega a^{\alpha \beta \sigma \tau} \gamma_{\sigma \tau}(\bm{\xi}) \gamma_{\alpha \beta}(\bm{\eta}) \sqrt{a} \dd y,\\
	B_F(\bm{\xi}, \bm{\eta})&:=\dfrac13 \int_\omega a^{\alpha \beta \sigma \tau} \rho_{\sigma \tau}(\bm{\xi}) \rho_{\alpha \beta}(\bm{\eta}) \sqrt{a} \dd y,
\end{align*}
for each $\bm{\xi}=(\xi_i)\in \bm{V}_K(\omega)$ and each $\bm{\eta}=(\eta_i) \in \bm{V}_K(\omega)$,
and define the linear form $\ell^\varepsilon$ by
$$
\ell^\varepsilon(\bm{\eta}):=\int_\omega p^{i,\varepsilon} \eta_i \sqrt{a} \dd y, \textup{ for each } \bm{\eta}=(\eta_i) \in \bm{V}_K(\omega), 
$$
where $p^{i,\varepsilon}(y):=\int_{-\varepsilon} ^\varepsilon f^{i,\varepsilon}(y,x_3) \dd x_3$ at each $y \in \omega$.

Then the \emph{total energy} of the shell is the \emph{quadratic functional} $J:\bm{V}_K(\omega) \to \mathbb{R}$ defined by
$$
J(\bm{\eta}):=\dfrac{\varepsilon}{2}B_M(\bm{\eta},\bm{\eta})+\dfrac{\varepsilon^3}{6}B_F(\bm{\eta},\bm{\eta})-\ell^\varepsilon(\bm{\eta})\quad\textup{ for each }\bm{\eta} \in \bm{V}_K(\omega).
$$

The terms $\displaystyle{\dfrac{\varepsilon}{2}B_M(\cdot,\cdot)}$ and $\displaystyle{\dfrac{\varepsilon^3}{6}B_F(\cdot,\cdot)}$ are respectively related to the \emph{membrane part} and the \emph{flexural part} of the total energy, as aptly recalled by the subscripts ``$M$'' and ``$F$''.

The unknown of Koiter's model is the ``two-dimensional'' displacement field $\zeta_{i,K}^\varepsilon \bm{a}^i:\overline{\omega} \to \mathbb{E}^3$ of the middle surface $\bm{\theta}(\overline{\omega})$ of the shell.
The vector field $\bm{\zeta}_K^\varepsilon:=(\zeta_{i,K}^\varepsilon)$ should thus be the solution of the following boundary value problem:

\begin{customprob}{$\mathcal{P}_K^\varepsilon(\omega)$}\label{Koiter}
	Find a vector field $\bm{\zeta}_K^\varepsilon:\overline{\omega} \to \mathbb{R}^3$ that satisfies
	$$
	\bm{\zeta}_K^\varepsilon \in \bm{V}_K(\omega) \quad\textup{ and }\quad J(\bm{\zeta}_K^\varepsilon)=\inf_{\bm{\eta} \in \bm{V}_K(\omega)} J(\bm{\eta}),
	$$
	\emph{or equivalently}, find $\bm{\zeta}_K^\varepsilon \in \bm{V}_K(\omega)$ that satisfies the following variational equations:
	$$
	\varepsilon B_M(\bm{\zeta}_K^\varepsilon, \bm{\eta})+\varepsilon^3 B_F(\bm{\zeta}_K^\varepsilon, \bm{\eta})=\ell^\varepsilon(\bm{\eta}) \quad\textup{ for all }\bm{\eta} \in \bm{V}_K(\omega).
	$$
	\bqed	
\end{customprob}

As first shown in~\cite{BerCia1976} (see also~\cite{BerCiaMia1994}), this problem has one and only one solution.

One of the most remarkable features of the model proposed by Koiter is that Koiter's equations are valid for all types of shells, even though it is clear that these equations cannot be recovered as the outcome of an asymptotic analysis of the three-dimensional equations of linearized elasticity. It is indeed noticeable that the thickness parameter $\varepsilon$ is no longer associated with the integration domain, but is now regarded as a multiplicative term for the \emph{membrane part} and the \emph{flexural part}. As we will see in the forthcoming sections, the powers of the parameter $\varepsilon$ entering Koiter's model, which were derived by Koiter by solely relying on assumption of mechanical and geometrical natures, will determine the properties the applied body forces and the applied surface forces will have to satisfy for allowing us to carry the asymptotic analysis out.

\section{Numerical study of linearly elastic elliptic membrane shells}\label{Sec:3}

Section \ref{Sec:2} was devoted to the formulation of the boundary value problem for a ``general'' linearly elastic shells. In what follows, we consider a specific class of shells, according to the following definition proposed in \cite{CiaLods1996b}.

A linearly elastic shell satisfying the various assumptions made in Section \ref{Sec:2} is said to be a \emph{linearly elastic elliptic membrane shell} if the following two additional assumptions are contemporarily satisfied: \emph{first}, $\gamma_0 = \gamma$, i.e., the homogeneous boundary condition of place holds along the \emph{entire lateral face} $\gamma \times \left[ - \varepsilon , \varepsilon \right]$ of the shell, and \emph{second}, the middle surface of the shell under consideration $\bm{\theta}(\overline{\omega})$ is \emph{elliptic} (cf. Section \ref{Sec:1}).

Consider the problem~\ref{problem0} for a family of linearly elastic elliptic membrane shells, all sharing the \emph{same middle surface} and whose thickness $2 \varepsilon > 0$ is considered as a ``small'' parameter approaching zero. 

In this case, it was shown by Ciarlet \& Lods~\cite{Ciarlet1996} that it is possible to recover a set of two-dimensional variational equations as a result of a rigorous asymptotic analysis. This methodology was first proposed by Ciarlet \& Destuynder in the seminal work~\cite{CiaDes1979}, and consists of the following steps.

First of all, each problem~\ref{problem0}, $\varepsilon > 0$ undergoes a scaling over a \emph{fixed domain} $\Omega$, which affects the \emph{unknowns} and \emph{assumptions on the data} that will have to be made. These scalings and assumptions depend on the type of shell under consideration; for example, the assumptions made for treating \emph{linearly elastic flexural shells} (cf., e.g., \cite{Ciarlet19963}) are different from the one presented in this section.

Contextually, we define the set
\[
\Omega := \omega \times \left] - 1, 1 \right[, \quad \Gamma_0 := \gamma_0 \times [-1,1], \quad \Gamma_{\pm}=\omega \times \{\pm 1\},
\]
where $x = (x_i)$ denotes a generic point in the set $\overline{\Omega}$, and $\partial_i := \partial/ \partial x_i$. Each point $x = (x_i)$ in $\overline{\Omega}$ corresponds to a unique point $x^\varepsilon = (x^\varepsilon_i)$ defined by
\[
x^\varepsilon_\alpha := x_\alpha = y_\alpha \text{ and } x^\varepsilon_3 := \varepsilon x_3.
\]
 By so doing, it turns out that $\partial^\varepsilon_\alpha = \partial_\alpha$ and $\partial^\varepsilon_3 = \displaystyle\frac{1}{\varepsilon} \partial_3$.

The unknown $\bm{u}^\varepsilon = (u^\varepsilon_i)$ that appears in the formulation of the variational problem~\ref{problem0} and denotes the displacement field a linearly elastic elliptic membrane shell undergoes under a deformation, is associated with the \emph{scaled unknown} $\bm{u} (\varepsilon) = (u_i(\varepsilon))$ by means of the following transformation
\[
u_i (\varepsilon) (x) := u^\varepsilon_i (x^\varepsilon)
\]
at each $x\in \overline{\Omega}$. 

Finally, we set the assumptions on applied body forces out: we \emph{assume} that there exist functions $f^i \in L^2(\Omega)$ and $h^i \in L^2(\Gamma_{+} \cup \Gamma_{-})$ \emph{independent of} $\varepsilon$ such that the following \emph{assumptions on the data} hold
\begin{equation}
\label{dataMembranes}
\begin{aligned}
f^{i, \varepsilon} (x^\varepsilon) &= f^i(x) \text{ at each } x \in \Omega,\\
h^{i, \varepsilon} (x^\varepsilon) &= \varepsilon h^i(x) \text{ at each } x \in \Gamma_{+} \cup \Gamma_{-}.
\end{aligned}
\end{equation}

It is clear that the independence of the Lam\'{e} constants on the thickness parameter $\varepsilon$ assumed in Section \ref{Sec:2} in the context of the  formulation of problem~\ref{problem0} implicitly constitutes another \emph{assumption on the data}.

As a result of the scalings, the following identities hold:
	\begin{align*}
	g(\varepsilon) (x) &:= g^\varepsilon (x^\varepsilon)
	\text{ and } A^{ijk\ell} (\varepsilon) (x) :=
	A^{ijk\ell, \varepsilon}(x^\varepsilon)
	\text{ at each } x\in \overline{\Omega}, \\
	e_{\alpha \| \beta} (\varepsilon; \bm{v}) &:= \frac12 (\partial_\beta v_\alpha + \partial_\alpha v_\beta ) - \Gamma^k_{\alpha \beta} (\varepsilon) v_k = e_{\beta \| \alpha} (\varepsilon; \bm{v}) , \\
	e_{\alpha \| 3} (\varepsilon; \bm{v}) = e_{3\| \alpha}(\varepsilon, \bm{v}) &:= \frac12 \left(\frac{1}{\varepsilon} \partial_3 v_\alpha + \partial_\alpha v_3 \right) - \Gamma^\sigma_{\alpha3} (\varepsilon) v_\sigma, \\
	e_{3\| 3} (\varepsilon ; \bm{v}) &:= \frac{1}{\varepsilon} \partial_3 v_3,
\end{align*}
where
\[
\Gamma^p_{ij} (\varepsilon) (x) := \Gamma^{p, \varepsilon}_{ij} (x^\varepsilon) \text{ at each } x\in \overline{\Omega}.
\]

We also let:
\[
\bm{g}_i (\varepsilon)(x) := g^\varepsilon_i (x^\varepsilon) \textup{ at each } x \in \Omega.
\]

The \emph{scaled} variational problem $\mathcal{P} (\varepsilon; \Omega)$ defined next constitutes the starting point of the asymptotic analysis performed in~\cite{Ciarlet1996}.

\begin{customprob}{$\mathcal{P} (\varepsilon; \Omega)$}\label{problem0s}
	Find $\bm{u}^\varepsilon \in \bm{V}(\Omega):= \{\bm{v} = (v_i) \in \bm{H}^1(\Omega) ; \; \bm{v} = \textbf{0} \text{ on } \Gamma_0\}$ that satisfies the following variational equations:
	$$
	\int_{\Omega}
	A^{ijk\ell}(\varepsilon) e_{k\| \ell} (\varepsilon;\bm{u}(\varepsilon)) e_{i\| j}(\varepsilon;\bm{v}) \sqrt{g(\varepsilon)} \dd x = \int_{\Omega} f^{i} v_i\sqrt{g(\varepsilon)} \dd x + \int_{\Gamma_{+} \cup \Gamma_{-}} h^{i} v_i \sqrt{g(\varepsilon)} \dd \Gamma,
	$$
	for all $\bm{v}^\varepsilon = (v^\varepsilon_i) \in \bm{V}(\Omega^\varepsilon)$.
	\bqed	
\end{customprob}

The variational problem~\ref{problem0s} is a re-writing of Problem~\ref{problem0} in terms of the newly defined scaled variables. Therefore, Problem~\ref{problem0s} admits a unique solution too.

The functions $e_{i\|j} (\varepsilon; \bm{v})$ appearing in Problem~\ref{problem0s} go by the name of \emph{scaled linearized strains in curvilinear coordinates}, and they correspond to the scaled displacement vector field $v_i \bm{g}^i(\varepsilon)$.

\subsection{Asymptotic analysis}

Consider a family of \emph{linearly elastic elliptic membrane shells} with thickness $2\varepsilon>0$ with each having the same middle surface $\bm{\theta}(\overline{\omega})$. It was shown by Ciarlet \& Lods~\cite{CiaLods1996b} that the solutions $\bm{u}(\varepsilon)=(u_i(\varepsilon))$ of the associated scaled three-dimensional problems~\ref{problem0s} converge in $H^1(\Omega) \times H^1(\Omega)\times L^2(\Omega)$ as $\varepsilon \to 0$ towards a \emph{limit} vector field $\bm{u}=(u_i)$. This limit can be proven to be independent of the transverse variable $x_3$, and can be identified with the solution $\overline{\bm{u}}$ of a two-dimensional variational problem posed over the set $\omega$.

The following result, due to Ciarlet \& Lods~\cite{Ciarlet1996}, is critical for carrying out the numerical experiments for linearly elastic elliptic membrane shells meant to validate the justification of Koiter's model for this very type of shells.

\begin{theorem}
\label{t:5}
Assume that $\bm{\theta} \in \mathcal{C}^2(\overline{\omega};\mathbb{R}^3)$. Consider a family of linearly elastic elliptic membrane shells with thickness $2\varepsilon$ approaching zero and with each having the same elliptic middle surface $\bm{\theta}(\overline{\omega})$, and let the assumptions on the data be as in Section~\ref{Sec:2} (for the Lam\'e constants) and~\eqref{dataMembranes}.

Let $\bm{u}(\varepsilon)$ denote, for $\varepsilon>0$ small enough, the solution of the associated
scaled three-dimensional problems~\ref{problem0s}. 

Let
\begin{align*}
	a^{\alpha\beta\sigma\tau}&=\dfrac{4\lambda\mu}{\lambda+2\mu} a^{\alpha\beta} a^{\sigma\tau}+2\mu (a^{\alpha\sigma} a^{\beta\tau} + a^{\alpha\tau} a^{\beta\sigma}),\\
	\gamma_{\alpha \beta}(\bm{\eta})&:=\frac12 (\partial_\beta \eta_\alpha + \partial_\alpha \eta_\beta ) - \Gamma^\sigma_{\alpha \beta} \eta_\sigma - b_{\alpha \beta} \eta_3,\\
	p^i&=\int_{-1}^{1} f^i \dd x_3 +h^i_{+}+h^i_{-}, \quad h^i_\pm=h^i(\cdot, \pm 1).
\end{align*}

Then there exist functions $u_\alpha \in H^1(\Omega)$ satisfying $u_\alpha=0$ on $\Gamma_0$, and a function $u_3 \in L^2(\Omega)$ such that:
\begin{align*}
u_\alpha(\varepsilon) \to u_\alpha\quad&\textup{ in } H^1(\Omega) \textup{ as }\varepsilon\to 0,\\
u_3(\varepsilon) \to u_3 \quad&\textup{ in } L^2(\Omega) \textup{ as }\varepsilon\to 0,\\
\bm{u}=(u_i) &\textup{ is independent of the transverse variable } x_3.
\end{align*}

Furthermore, the \emph{average}
$$
\overline{\bm{u}}:=\dfrac{1}{2}\int_{-1}^{1} \bm{u} \dd x_3
$$
satisfies the following scaled variational two-dimensional problem $\mathcal{P}_M(\omega)$ of a linearly elastic elliptic membrane shell:
\begin{customprob}{$\mathcal{P}_M(\omega)$}\label{limitM}
	Find $\overline{\bm{u}} \in \bm{V}_M(\omega):= H^1_0(\omega) \times H^1_0(\omega) \times L^2(\omega)$ that satisfies the following variational equations
	$$
	\int_{\omega}
	a^{\alpha\beta\sigma\tau} \gamma_{\sigma\tau}(\overline{\bm{u}}) \gamma_{\alpha \beta}(\bm{\eta}) \sqrt{a} \dd y = \int_{\omega} p^{i} \eta_i\sqrt{a} \dd y,
	$$
	for all $\bm{\eta} = (\eta_i) \in \bm{V}_M(\omega)$.
	\bqed	
\end{customprob}
and is the unique vector field doing so.
\qed
\end{theorem}

Once this asymptotic analysis is carried out, the next step consists in ``\emph{de-scaling}'' the results of Theorem~\ref{t:5}, which apply to the solutions $\bm{u}(\varepsilon)$ of the \emph{scaled} problem~\ref{problem0s}.

This means that these results have to be read in terms of the \emph{unknown} $u^\varepsilon_i \bm{g}^{i, \varepsilon} : \overline{\Omega^\varepsilon} \to \mathbb{E}^3$, which represents the physical \emph{three-dimensional vector field} of the actual reference configuration of the shell, viz., the one having thickness equal to $2\varepsilon$. 

The next result shows that this operation is smoothly performed through the exploitation of the special \emph{averages} $\displaystyle \frac{1}{2\varepsilon} \int^\varepsilon_{-\varepsilon} u^\varepsilon_i \bm{g}^{i, \varepsilon} \dd x^\varepsilon_3$ across the thickness of the shell.

\begin{theorem}[Theorem 4.6-1 in \cite{Ciarlet2000}] \label{t:6}
	Let the assumptions on the data be as in \emph{Section \ref{Sec:3}} and let $\bm{\theta} \in \mathcal{C}^3 (\overline{\omega}; \mathbb{E}^3)$ be a smooth enough immersion. Let $\bm{u}^\varepsilon = (u^\varepsilon_i) \in \bm{V}(\Omega^\varepsilon)$ denote for each $\varepsilon > 0$ the unique solution of the variational problem~\ref{problem0}, and let $\bm{\zeta} = (\zeta_i) \in \bm{V}_M(\omega)$ denote the unique solution for Problem~\ref{limitM}. Then
	\begin{align*}
		\frac{1}{2\varepsilon} \int^\varepsilon_{-\varepsilon} u^\varepsilon_\alpha \bm{g}^{\alpha , \varepsilon} \dd x^\varepsilon_3 & \to \zeta_\alpha \bm{a}^\alpha \textup{ in } \bm{H}^1 (\omega) \textup{ as } \varepsilon \to 0, \\
		\frac{1}{2\varepsilon} \int^\varepsilon_{-\varepsilon} u^\varepsilon_3 \bm{g}^{3 , \varepsilon} \dd x^\varepsilon_3 &\to \zeta_3 \bm{a}^3
		\textup{ in } \bm{L}^2 (\omega) \textup{ as } \varepsilon \to 0. 
	\end{align*}
\qed
\end{theorem}

\subsection{Justification of Koiter's model for linearly elastic elliptic membrane shells}
\label{S:1:1}

The justification of Koiter's model for linearly elastic elliptic membrane shells was established by Ciarlet \& Lods in~\cite[Theorem~2.1]{Ciarlet19962}.

\begin{theorem}
\label{t:k:1}
Assume that $\bm{\theta} \in \mathcal{C}^2(\overline{\omega};\mathbb{R}^3)$. Consider a family of linearly elastic elliptic membrane shells with thickness $2\varepsilon$ approaching zero and with each having the same elliptic middle surface $\bm{\theta}(\overline{\omega})$, and let the assumptions on the data be as in Section~\ref{Sec:2} (for the Lam\'e constants) and~\eqref{dataMembranes}. For each $\varepsilon>0$ let $\bm{u}^\varepsilon=(u_i^\varepsilon) \in \bm{H}^1(\Omega^\varepsilon)$ and $\bm{\zeta}_K^\varepsilon=(\zeta_{i,K}^\varepsilon)$ denote the solutions to the three-dimensional problem~\ref{problem0} and the two-dimensional problem~\ref{Koiter}.
Let also
$$
\bm{\zeta} =(\zeta_i) \in \bm{V}_M(\omega)= H^1_0(\omega) \times H^1_0(\omega) \times L^2(\omega),
$$
denote the solution to the two-dimensional scaled variational problem~\ref{limitM}, solution which is independent of $\varepsilon$. Then, the following convergences hold:
\begin{align*}
\dfrac{1}{2\varepsilon} \int_{-\varepsilon}^{\varepsilon} u^\varepsilon_\alpha \bm{g}^{\alpha,\varepsilon} \dd x_3^\varepsilon \to \zeta_\alpha \bm{a}^\alpha &\textup{ in } \bm{H}^1(\omega) \textup{ as } \varepsilon \to 0,\\
\zeta_{\alpha,K}^\varepsilon \bm{a}^\alpha \to \zeta_\alpha \bm{a}^\alpha &\textup{ in } \bm{H}^1(\omega) \textup{ as } \varepsilon \to 0,
\end{align*}
and
\begin{align*}
	\dfrac{1}{2\varepsilon} \int_{-\varepsilon}^{\varepsilon} u^\varepsilon_3 \bm{g}^{3,\varepsilon} \dd x_3^\varepsilon \to \zeta_3 \bm{a}^3 &\textup{ in } \bm{L}^2(\omega) \textup{ as } \varepsilon \to 0,\\
	\zeta_{3,K}^\varepsilon \bm{a}^3 \to \zeta_3 \bm{a}^3 &\textup{ in } \bm{L}^2(\omega) \textup{ as } \varepsilon \to 0.
\end{align*}
\qed
\end{theorem}

Note that the first and the third convergence in Theorem~\ref{t:k:1} have been recovered in Theorem~\ref{t:6}.

\subsection{Numerical results}
\label{S:1:2}

We verify the convergences announced in Theorem~\ref{t:6} through a series of \emph{ad hoc} numerical experiments carried out by implementing the finite element method (cf., e.g., \cite{PGCFEM}) via the software FreeFem~\cite{Hecht2012}. The results are visualised in ParaView~\cite{Ahrens2005}. Specifically, we use conforming finite element to discretize the three components of the displacement(cf., e.g., \cite{Shen2021}). Apart from the transverse component of the solution of Koiter's model, which is discretized by means of HCT triangles (cf., e.g., Chapter~6 of~\cite{PGCFEM}), all the other components of the solution of Koiter's model, the solution of the three-dimensional model and the solution of the two-dimensional limit model are approximated via a Lagrange finite element.

We use a portion of a spherical shell for numerical experiments and assume that the thickness of the shell is $2\varepsilon$. 

The domain $\omega$ is defined as follows
\begin{align*}
	\omega:=\left\{(y_1,y_2)\in\mathbb{R}^2;\left(\frac{\pi}{6},\frac{5\pi}{6}\right)\times(0, \pi)\right\},
\end{align*}
where $m=0.06\meter,n=0.05\meter,l=0.03\meter$, and $\gamma_0$ is the region of the boundary at which the clamping occur.

According to \cite{Shen2021}, in curvilinear coordinates, the middle surface is given by the mapping $\bm{\theta}$ defined by
\begin{align*}
	\bm{\theta}(y_1,y_2)=(m\sin y_1 \cos y_2,n \sin y_1 \sin y_2,l \cos y_1),\quad\textup{ for all } (y_1,y_2)\in\overline{\omega}.
\end{align*}

The first experiment we carry out concerns the displacement of a three-dimensional linearly elastic elliptic membrane shell whose middle surface is a half-sphere. We use the following values for the Lam\'e constants , Young's modulus, the Poission ratio and the applied body force density (cf., e.g., \cite{Shen2021}):
\begin{align*}
\nu&=0.25,\\
E&=2.0\times10^{11}\pascal,\\
\lambda&=8.0\times10^{10}\pascal,\\
\mu&=8.0\times10^{10}\pascal,\\
f^{i,\varepsilon}(x^\varepsilon)&=1.0\times10^{-1}\newton,\\
h^{i,\varepsilon}(x^\varepsilon)&=0\newton.
\end{align*}

In Table \ref{table:va1}, ErrLK indicates the residual error between the solution of the two-dimensional(2D) limit model and the solution of Koiter's model; Err3DL indicates the residual error between the solution of the 2D limit model and the solution of the three-dimensional(3D) linearly elastic elliptic membrane shell model; Err3DK indicates the residual error between the solution of Koiter's model and the solution of the 3D linearly elastic elliptic membrane shell model. The consistent change in Table \ref{table:va1} means that as the thickness parameter $\varepsilon$ approaches zero, the residual errors between different models gradually tend to zero as well.

\begin{figure}[H]
    \centering
    \begin{subfigure}[b]{0.33\linewidth}
		\includegraphics[width=1.0\linewidth]{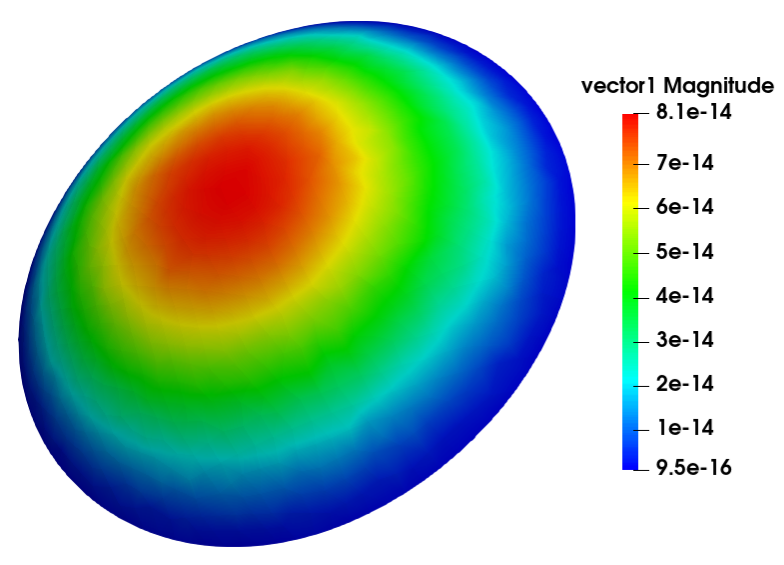}
		\subcaption{$\varepsilon=8\times10^{-2}$}
	\end{subfigure}%
	\hspace{0.5cm}
	\begin{subfigure}[b]{0.33\linewidth}	
	    \includegraphics[width=1.0\linewidth]{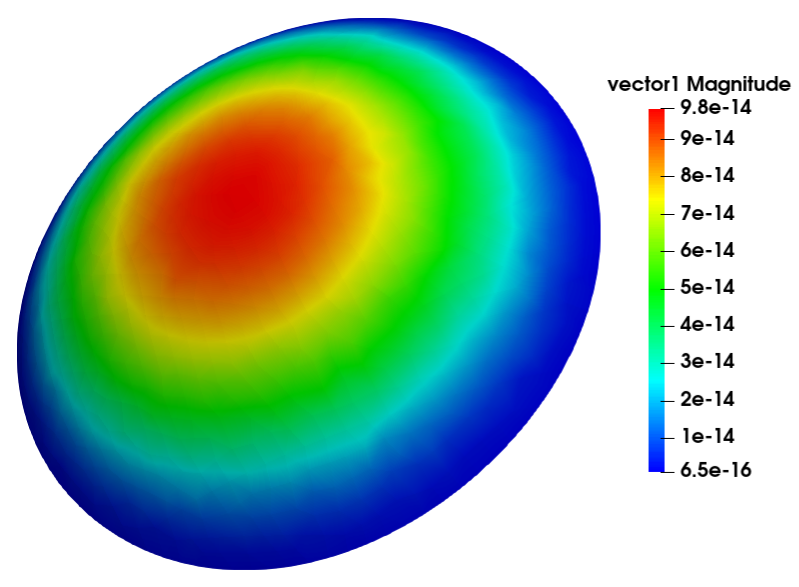}
		\subcaption{$\varepsilon=4\times10^{-2}$}
	\end{subfigure}%
	\vspace{0.5cm}
	\begin{subfigure}[b]{0.33\linewidth}
		\includegraphics[width=1.0\linewidth]{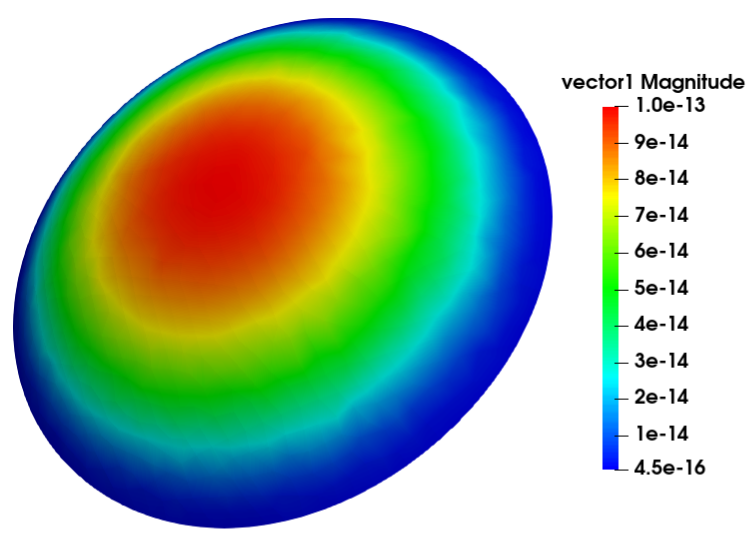}
		\subcaption{$\varepsilon=2\times10^{-2}$}
	\end{subfigure}%
	\hspace{0.5cm}
	\begin{subfigure}[b]{0.33\linewidth}
		\includegraphics[width=1.0\linewidth]{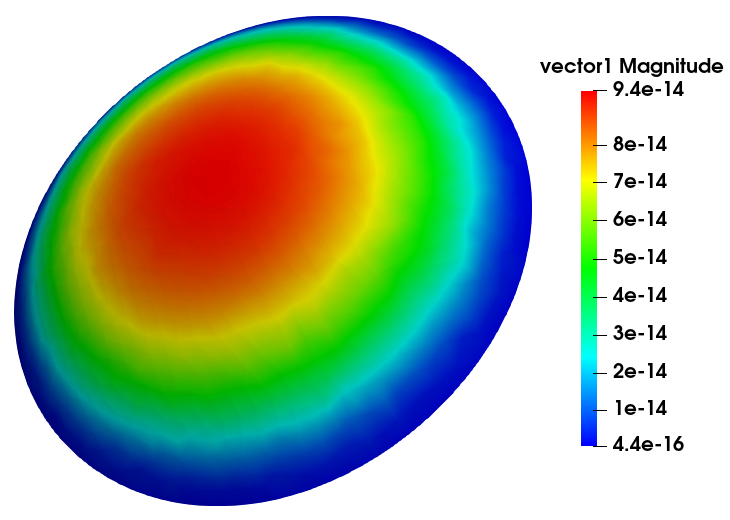}
		\subcaption{$\varepsilon=1\times10^{-2}$}
	\end{subfigure}%
	\caption{Displacement of a three-dimensional linearly elastic elliptic membrane shell whose middle surface is a half-sphere. Different values of the parameter $\varepsilon$ are taken into account.}
	\label{EMS1}
\end{figure}

\begin{figure}[H]
	\centering
	\begin{subfigure}[b]{0.33\linewidth}
		\includegraphics[width=1.0\linewidth]{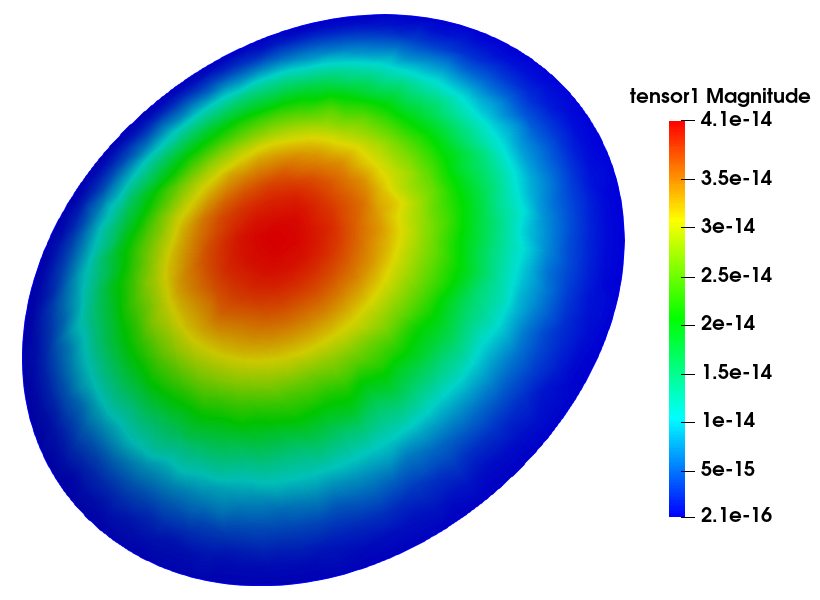}
		\subcaption{$\varepsilon=8\times10^{-2}$}
	\end{subfigure}%
	\hspace{0.5cm}
	\begin{subfigure}[b]{0.33\linewidth}
		\includegraphics[width=1.0\linewidth]{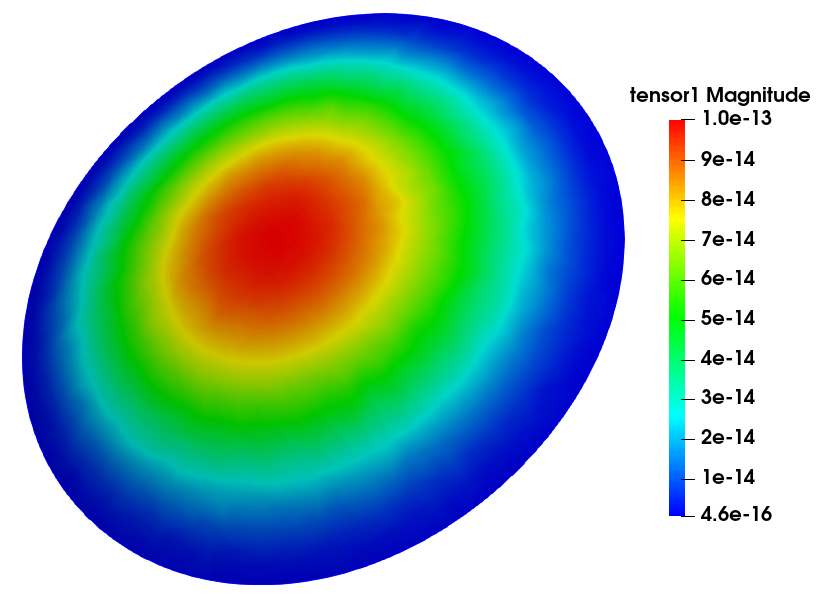}
		\subcaption{$\varepsilon=4\times10^{-2}$}
	\end{subfigure}%
	\vspace{0.5cm}
	\begin{subfigure}[b]{0.33\linewidth}
		\includegraphics[width=1.0\linewidth]{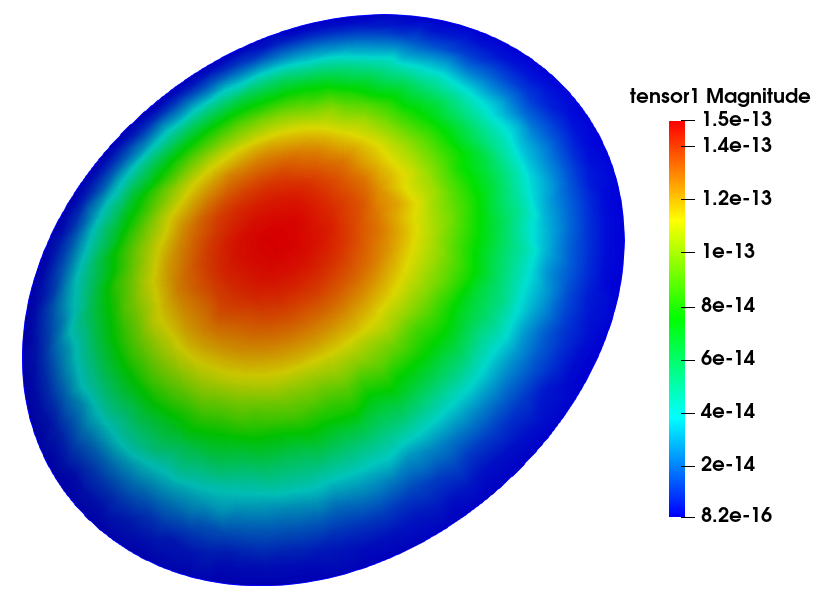}
		\subcaption{$\varepsilon=2\times10^{-2}$}
	\end{subfigure}%
    \hspace{0.5cm}
	\begin{subfigure}[b]{0.33\linewidth}
		\includegraphics[width=1.0\linewidth]{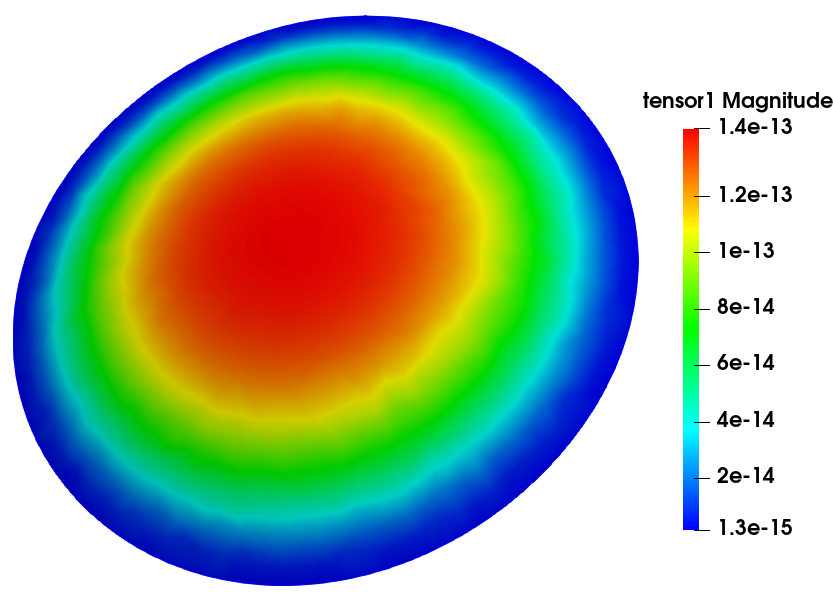}
		\subcaption{$\varepsilon=1\times10^{-2}$}
	\end{subfigure}%
	\caption{Displacement of a two-dimensional Koiter's shell under the assumptions on the data~\eqref{dataMembranes} characterizing linearly elastic elliptic membrane shells. Different values of the thickness parameter $\varepsilon$ are taken into account.}
	\label{EMS2}
\end{figure}

Figure \ref{EMS1} and Figure \ref{EMS2} show the numerical results of the 3D linearly elastic elliptic membrane shell model and Koiter's model. The entire bottom of each figure is in blue color, indicating that the entire boundary is clamped, whereas the biggest deformation happens on the parts in red color, i.e., the center of the ellipsoidal shell. It can be observed from the consistent shapes of the figures that the 3D linearly elastic elliptic membrane shell model and Koiter's model manifest the same asymptotic behavior as the thickness parameter $\varepsilon$ approaches zero.

\begin{table}[htbp]
\centering
\caption{Spherical elliptic membrane shell result.}
\hspace{0.5cm}
\Huge{
 \resizebox{\textwidth}{19mm}{\setlength{\tabcolsep}{22mm}{
  \begin{tabular}{r@{\,}l c c c}
  \toprule[ 2pt]
  \specialrule{0em}{3pt}{3pt}
   \multicolumn{2}{c}{$\varepsilon$}&\multicolumn{1}{c}{ErrLK}&\multicolumn{1}{c}{Err3DK}&\multicolumn{1}{c}{Err3DL}\\
   \midrule[1pt]
    8&$\times10^{-2}$    & 1.14913e$-$14 & 8.47566e$-$15 & 1.9967e$-$14    \\
    4&$\times10^{-2}$    & 5.5481e$-$15  & 8.38406e$-$15 & 1.39322e$-$14   \\
    2&$\times10^{-2}$    & 2.21322e$-$15 & 7.58495e$-$15 & 9.79817e$-$15   \\
    1&$\times10^{-2}$    & 9.06036e$-$16 & 5.95969e$-$15 & 6.86572e$-$15   \\
    5&$\times10^{-3}$    & 3.90062e$-$16 & 4.42482e$-$15 & 4.81488e$-$15   \\
    2.5&$\times10^{-3}$  & 1.67942e$-$16 & 3.21773e$-$15 & 3.38566e$-$15   \\
    1.25&$\times10^{-3}$ & 7.08079e$-$17 & 2.3155e$-$15  & 2.38631e$-$15   \\
    6.25&$\times10^{-4}$ & 2.851e$-$17   & 1.65595e$-$15 & 1.68446e$-$15   \\
  \bottomrule[ 2pt]
  \end{tabular}\label{table:va1}
}}}
\end{table}

\section{Numerical study of generalized membrane shells of the first kind}\label{Sec:4}

We now consider a linearly elastic shell subjected to the various assumptions announced in Section \ref{Sec:2}. One such shell is said to be a \emph{linearly elastic generalized membrane shell} (from now onward simply \emph{generalized membrane shell}) if the following three additional assumptions are satisfied: \emph{first}, the shell is subjected to a boundary condition of place along
a portion of its lateral face with $\bm{\theta}(\gamma_0)$ as its middle curve, where the subset $\gamma_0 \subset \gamma$ has nonzero length, \emph{second}, the space of linearised inextensional displacements
$$
\bm{V}_F(\omega):=\{\bm{\eta}=(\eta_i) \in H^1(\omega) \times H^1(\omega)\times H^2(\omega); \eta_i=\partial_{\nu}\eta_3=0 \textup{ on }\gamma_0\},
$$
is only made of the null vector field, i.e., 
$$
\bm{V}_F(\omega)=\{\textbf{0}\},
$$
and, \emph{third}, the shell is \emph{not} a linearly elastic elliptic membrane shells, in the sense that \emph{either} $\gamma_0 \subsetneq \gamma$ \emph{or} the middle surface $\bm{\theta}(\overline{\omega})$ is not elliptic.

Generalized membrane shells thus exhaust all the remaining cases of linearly elastic membrane shells, i.e., those for which $\bm{V}_F(\omega)=\{\textbf{0}\}$.

As in section~\ref{Sec:3}, wthe first step consists in scaling each problem~\ref{problem0}, $\varepsilon > 0$, over a \emph{fixed domain} $\Omega$, by making use of appropriate \emph{scalings on the unknowns} and \emph{assumptions on the data}.

In this direction, we define the set
\[
\Omega := \omega \times \left] - 1, 1 \right[, \quad \Gamma_0 := \gamma_0 \times [-1,1], \quad \Gamma_{\pm}=\omega \times \{\pm 1\},
\]
and we let $x = (x_i)$ denote a generic point in the set $\overline{\Omega}$, and $\partial_i := \partial/ \partial x_i$. Each point $x = (x_i) \in \overline{\Omega}$ is associated with the uniquely determined point $x^\varepsilon = (x^\varepsilon_i)$ defined by
\[
x^\varepsilon_\alpha := x_\alpha = y_\alpha \text{ and } x^\varepsilon_3 := \varepsilon x_3,
\]
so that $\partial^\varepsilon_\alpha = \partial_\alpha$ and $\partial^\varepsilon_3 = \displaystyle\frac{1}{\varepsilon} \partial_3$. 

The unknown $\bm{u}^\varepsilon = (u^\varepsilon_i)$ and the vector fields $\bm{v}^\varepsilon = (v^\varepsilon_i)$ entering the formulation of the problem~\ref{problem0} are associated, thanks to the scaling announced beforehand, with the scaled unknown $\bm{u} (\varepsilon) = (u_i(\varepsilon))$ and the scaled vector fields $\bm{v} = (v_i)$ defined by
\[
u_i (\varepsilon) (x) := u^\varepsilon_i (x^\varepsilon) \text{ and } v_i(x) := v^\varepsilon_i (x^\varepsilon)
\]
at each $x\in \overline{\Omega}$.

As a result of the announced scalings, the following identities hold (cf. section~\ref{Sec:3}):
\begin{align*}
	g(\varepsilon) (x) &= g^\varepsilon (x^\varepsilon)
	\text{ and } A^{ijk\ell} (\varepsilon) (x) = A^{ijk\ell, \varepsilon}(x^\varepsilon) \text{ at each } x\in \overline{\Omega}, \\
	e_{\alpha \| \beta} (\varepsilon; \bm{v}) &= \frac12 (\partial_\beta v_\alpha + \partial_\alpha v_\beta ) - \Gamma^k_{\alpha \beta} (\varepsilon) v_k = e_{\beta \| \alpha} (\varepsilon; \bm{v}) , \\
	e_{\alpha \| 3} (\varepsilon; \bm{v}) = e_{3\| \alpha}(\varepsilon, \bm{v}) &= \frac12 \left(\frac{1}{\varepsilon} \partial_3 v_\alpha + \partial_\alpha v_3 \right) - \Gamma^\sigma_{\alpha3} (\varepsilon) v_\sigma, \\
	e_{3\| 3} (\varepsilon ; \bm{v}) := \frac{1}{\varepsilon} \partial_3 v_3,
\end{align*}
where
\[
\Gamma^p_{ij} (\varepsilon) (x) = \Gamma^{p, \varepsilon}_{ij} (x^\varepsilon) \text{ at each } x\in \overline{\Omega}.
\]

We also define:
\[
\bm{g}_i (\varepsilon)(x) = g^\varepsilon_i (x^\varepsilon) \textup{ at each } x \in \Omega.
\]

The \emph{second} condition appearing in the definition of a generalized membrane shell, namely,
$$
\bm{V}_F(\omega)=\{\bf0\},
$$
is equivalent to stating that the \emph{semi-norm $\left|\cdot\right|_\omega^M$ defined by
	$$
	|\bm{\eta}|_\omega^M:=\bigg\{\sum_{\alpha, \beta}|\gamma_{\alpha\beta}(\bm{\eta})|_{0,\omega}^2\bigg\}^{1/2},
	$$
	for each $\bm{\eta}=(\eta_i) \in H^1(\omega) \times H^1(\omega) \times L^2(\omega)$ is actually a norm over the space} 
$$
\bm{V}_K(\omega)=\{\bm{\eta}=(\eta_i) \in H^1(\omega)\times H^1(\omega) \times H^2(\omega); \eta_i=\partial_{\nu}\eta_3=0 \textup{ on }\gamma_{0}\}.
$$ 

Generalized membranes are further classified into two sub-categories by means of the spaces:
\begin{align*}
	\bm{V}(\omega)&:=\{\bm{\eta}=(\eta_i) \in \bm{H}^1(\omega);\bm{\eta}={\bf0} \textup{ on }\gamma_{0}\},\\
	\bm{V}_{0}(\omega)&:=\{\bm{\eta}=(\eta_i)\in\bm{V}(\omega);\gamma_{\alpha \beta}(\bm{\eta})=0 \textup{ in }\omega\}.
\end{align*}
A generalized membrane shell is \emph{``of the first kind''} if 
$$
\bm{V}_{0}(\omega)=\{\bf0\},
$$
or, equivalently, \emph{if the semi-norm} $\left|\cdot\right|_\omega^M$ \emph{is already a norm over the space} $\bm{V}(\omega)$ (hence, \emph{a fortiori}, over the space $\bm{V}_K(\omega) \subset \bm{V}(\omega)$).

Otherwise, namely, if 
$$
{\bm{V}_F(\omega)=\{\bf0\}} \quad\textup{but}\quad {\bm{V}_{0}(\omega)\neq \{\bf0\}},
$$
or, equivalently, if the semi-norm $\left|\cdot\right|_\omega^M$ actually becomes a norm over the space $\bm{V}_K(\omega)$ but \emph{not} over the space $\bm{V}(\omega)$, the linearly elastic shell is said to be a generalized membrane shell ``\emph{of the second kind}".
In this paper we shall only consider generalized membranes of the ``first kind" (which are the most frequently encountered in practice).

The space
$$
\bm{V}_0(\Omega):=\{\bm{v}\in \bm{H}^1(\Omega); \bm{v}=\bm{0} \textup{ on }\Gamma_{0}, \partial_3 \bm{v}=\bm{0} \textup{ in }\Omega \textup{ and }\gamma_{\alpha \beta}(\overline{\bm{v}})=0 \textup{ in }\omega\},
$$
where, as usual, $\overline{\bm{v}}$ denotes the average of $\bm{v}$ with respect to the transverse variable, will play a crucial role in the recovery of the set of two-dimensional equations via the rigorous asymptotic analysis announced later on in this section. We observe that the space $\bm{V}_0(\Omega)$ is the ``three-dimensional analogue'' of the space $\bm{V}_0(\omega)$. 

In a similar fashion, define the semi-norm $|\cdot|_\Omega^M$ by
$$
|\bm{v}|_\Omega^M:=\{\|\partial_3 \bm{v}\|_{\bm{L}^2(\Omega)}^2 + (|\overline{\bm{v}}|_\omega^M)^2\}^{1/2},\quad\textup{ for all }\bm{v}\in \bm{V}(\Omega),
$$
which is the ``three-dimensional analogue'' of the semi-norm $|\cdot|_\omega^M$. 

It is also possible to observe that
$$
\bm{V}_0(\Omega)=\{\bm{0}\} \textup{ if and only if } \bm{V}_0(\omega)=\{\bm{0}\},
$$
or, equivalently, that $|\cdot|_\omega^M$ is a norm in the space $\bm{V}(\omega)$ if and only if $|\cdot|_\Omega^M$ is a norm in the space $\bm{V}(\Omega)$.

To sum up, a linearly elastic generalized membrane shell is of the first kind if and only if $\bm{V}_0(\Omega)=\{\bm{0}\}$.

In order to conduct the asymptotic analysis for a family of generalized membrane shells of the first kind, we need to assume
that the applied forces are \emph{admissible}, in the sense that there exists a constant $\kappa_0>0$ \emph{independent of $\varepsilon$} such that:
$$
|L(\varepsilon)(\bm{v})| \le \kappa_0 \left\{\sum_{i,j}\|e_{i\|j}(\varepsilon;\bm{v})\|_{L^2(\Omega)}\right\}^{1/2},\quad\textup{ for all } 0<\varepsilon\le \varepsilon_{1} \textup{ and all } \bm{v} \in \bm{V}(\Omega),
$$
where $L(\varepsilon)$ is the linear form appearing in the variational formulation of the corresponding three-dimensional model. The details concerning the admissibility of applied body forces are discussed in Chapter~5 of~\cite{Ciarlet2000}.

\subsection{Asymptotic analysis}

We are now in a position to state the convergence theorem for a family of generalized membrane shells of the first kind, that was originally proposed in~\cite{Ciarlet19964}, and was later on improved in~\cite[Theorem~5.6-1]{Ciarlet2000}.

\begin{theorem}
\label{t:8}
Assume that $\bm{\theta} \in \mathcal{C}^3(\overline{\omega};\mathbb{R}^3)$. Consider a family of generalized membrane shells of the first kind with thickness $2\varepsilon$ approaching zero, with each having the same middle surface $\bm{\theta}(\overline{\omega})$, with each subjected to a boundary condition of place along a portion of its lateral face having the same set $\bm{\theta}(\gamma_0)$ as its middle curve, and subjected to applied forces that are admissible.
Let $\bm{u}(\varepsilon)$ denote, for each $\varepsilon>0$ small enough, the solution of the associated scaled three-dimensional problems~\ref{problem0s}. Define the spaces:
\begin{align*}
\bm{V}_M^\sharp(\Omega):= \textup{ abstract completion of } \bm{V}(\Omega) \textup{ with respect to } |\cdot|_\Omega^M,\\
\bm{V}_M^\sharp(\omega):= \textup{ abstract completion of } \bm{V}(\omega) \textup{ with respect to } |\cdot|_\omega^M.
\end{align*}

Then, there exists $\bm{u} \in \bm{V}_M^\sharp(\Omega)$ and there exists $\bm{\zeta} \in \bm{V}_M^\sharp(\omega)$ such that
\begin{align*}
\bm{u}(\varepsilon) \to \bm{u} \textup{ in }\bm{V}_M^\sharp(\Omega) &\textup{ as } \varepsilon \to 0,\\
\overline{\bm{u}}(\varepsilon):=\dfrac{1}{2}\int_{-1}^{1} \bm{u}(\varepsilon) \dd x_3 \to \bm{\zeta} \textup{ in }\bm{V}_M^\sharp(\omega) &\textup{ as } \varepsilon \to 0.
\end{align*}

Let the assumptions on the data be as in Section~\ref{Sec:2} (for the Lam\'e constants) and let
\begin{align*}
a^{\alpha\beta\sigma\tau}&:=\dfrac{4\lambda\mu}{\lambda+2\mu} a^{\alpha\beta} a^{\sigma\tau}+2\mu (a^{\alpha\sigma} a^{\beta\tau} + a^{\alpha\tau} a^{\beta\sigma}),\\
\gamma_{\alpha \beta}(\bm{\eta})&:=\frac12 (\partial_\beta \eta_\alpha + \partial_\alpha \eta_\beta ) - \Gamma^\sigma_{\alpha \beta} \eta_\sigma - b_{\alpha \beta} \eta_3,\\
B_M(\bm{\xi},\bm{\eta})&:=\int_{\omega} a^{\alpha\beta\sigma\tau} \gamma_{\sigma\tau}(\bm{\xi}) \gamma_{\alpha \beta}(\bm{\eta}) \sqrt{a} \dd y, \textup{ for all }\bm{\xi}, \bm{\eta} \in \bm{V}(\omega),\\
L_M(\bm{\eta})&:=\int_{\omega} \varphi^{\alpha\beta} \gamma_{\alpha \beta}(\bm{\eta}) \sqrt{a} \dd y,\\
\varphi^{\alpha\beta}&:=\int_{-1}^{1}\left\{F^{\alpha\beta}-\dfrac{\lambda}{\lambda+2\mu}a^{\alpha\beta}F^{33}\right\} \dd x_3 \in L^2(\omega),
\end{align*}
where the functions $F^{ij} \in L^2(\Omega)$ are those used in the definition of admissible forces, and $B_M^\sharp$ and $L_M^\sharp$ denote the unique continuous extensions from $\bm{V}(\omega)$ to $\bm{V}_M^\sharp(\omega)$ of the bilinear form $B_M$ and the linear form $L_M$. Then the limit $\bm{\zeta}$ satisfies the following scaled two-dimensional variational problem $\mathcal{P}_M^\sharp(\omega)$ of a generalized membrane shell of the first kind:
\begin{customprob}{$\mathcal{P}_M^{\sharp}(\omega)$}\label{limitGM}
	Find $\bm{\zeta}\in \bm{V}_M^\sharp(\omega)$ that satisfies the following variational equations:
	$$
	B_M^\sharp(\bm{\zeta},\bm{\eta})=L_M^\sharp(\bm{\eta}),
	$$
	for all $\bm{\eta} = (\eta_i) \in \bm{V}_M^\sharp(\omega)$.
	\bqed	
\end{customprob}
The vector field $\bm{\zeta}$ is the one and only one which solves Problem~\ref{limitGM}.
\qed
\end{theorem}

Like in section~\ref{Sec:3}, we ``\emph{de-scale}'' the results of Theorem \ref{t:8}, which apply to the solutions $\bm{u}(\varepsilon)$ of the \emph{scaled} problem~\ref{problem0s}. This means that we need to translate these results into ones about the \emph{unknown} $u^\varepsilon_i \bm{g}^{i, \varepsilon} : \overline{\Omega^\varepsilon} \to \mathbb{E}^3$, which represents the physical \emph{three-dimensional vector field} of the actual reference configuration of the shell. As shown in the next theorem, this operation is performed through the introduction of the \emph{averages} $\displaystyle \frac{1}{2\varepsilon} \int^\varepsilon_{-\varepsilon} u^\varepsilon_i \bm{g}^{i, \varepsilon} \dd x^\varepsilon_3$ across the thickness of the shell.

\begin{theorem}[Theorem 5.8-1 in \cite{Ciarlet2000}] \label{t:9}
	Let the assumptions on the data be as in \emph{Section \ref{Sec:3}} and let the assumptions on the immersion $\bm{\theta} \in \mathcal{C}^3 (\overline{\omega}; \mathbb{E}^3)$ be as in \emph{Theorem \ref{t:8}}. Let $\bm{u}^\varepsilon = (u^\varepsilon_i) \in \bm{V}(\Omega^\varepsilon)$ denote for each $\varepsilon > 0$ the unique solution of the variational problem~\ref{problem0}, and let $\bm{\zeta} = (\zeta_i) \in \bm{V}_M(\omega)$ denote the unique solution for Problem~\ref{limitGM}. Then
	\begin{align*}
		\frac{1}{2\varepsilon} \int^\varepsilon_{-\varepsilon} \bm{u}^\varepsilon \dd x^\varepsilon_3 & \to \bm{\zeta} \textup{ in } \bm{V}_M^\sharp(\omega) \textup{ as } \varepsilon \to 0.
	\end{align*}
	\qed
\end{theorem}

\subsection{Justification of Koiter's model for generalized membrane shells of the first kind}
\label{S:2:1}
The justification of Koiter's model for generalized membrane shells was established by Ciarlet \& Lods in the paper~\cite[Theorems~6.1 and~6.2]{Ciarlet19964}.
The preamble is the same as in section~\ref{S:1:1}.

\begin{theorem}
	\label{t:k:2}
	Assume that $\bm{\theta} \in \mathcal{C}^2(\overline{\omega};\mathbb{R}^3)$. Consider a family of generalized membrane shells with thickness $2\varepsilon$ approaching zero and with each having the same elliptic middle surface $\bm{\theta}(\overline{\omega})$, and let the assumptions on the data be as in Section~\ref{Sec:2}. For each $\varepsilon>0$ let $\bm{u}^\varepsilon=(u_i^\varepsilon) \in \bm{H}^1(\Omega^\varepsilon)$ and $\bm{\zeta}_K^\varepsilon=(\zeta_{i,K}^\varepsilon)$ denote the solutions to the three-dimensional problem~\ref{problem0} and the two-dimensional problem~\ref{Koiter}.
	Let also
	\begin{align*}
	\bm{\zeta} =(\zeta_i) \in \bm{V}_M^\sharp(\omega)&:=\textup{ abstract completion of the space }\bm{V}(\omega) \textup{ with repsect to } |\cdot|_\omega^M,\\
	\bm{V}(\omega)&:=\{\bm{\eta}=(\eta_i) \in \bm{H}^1(\omega);\bm{\eta}=\bm{0} \textup{ on }\gamma_{0}\},\\
	|\bm{\eta}|_\omega^M&:=\left\{\sum_{\alpha, \beta} \|\gamma_{\alpha \beta}(\bm{\eta})\|_{L^2(\omega)}^2\right\}^{1/2},
	\end{align*}
	denote the solution to the two-dimensional scaled variational problem~\ref{limitGM}, solution which is independent of $\varepsilon$. Then, the following convergences hold:
	\begin{align*}
		\dfrac{1}{2\varepsilon} \int_{-\varepsilon}^{\varepsilon} \bm{u}^\varepsilon \dd x_3^\varepsilon \to \zeta_\alpha \bm{a}^\alpha &\textup{ in } \bm{V}_M^\sharp(\omega) \textup{ as } \varepsilon \to 0,\\
		\bm{\zeta}_{K}^\varepsilon \to \bm{\zeta} &\textup{ in } \bm{V}_M^\sharp(\omega) \textup{ as } \varepsilon \to 0.
	\end{align*}
	\qed
\end{theorem}

Note that the first and the third convergence in Theorem~\ref{t:k:2} have been recovered in Theorem~\ref{t:9}.

\subsection{Numerical results}
Same as above, we can verify the convergences announced in Theorem~\ref{t:9} through a series of \emph{ad hoc} numerical experiments carried out by implementing the finite element method. Unlike linearly elastic elliptic membrane shells (cf. section~\ref{Sec:3}), we observe that the two-dimensional limit model is posed over an abstract completion of a Sobolev space. In general, this space is quite hard to exhibit and is not, \emph{a priori}, a subspace of the space of distributions. It was proved by Mardare~\cite{Mardare1998} that this space is a subspace of the space of distributions under special geometrical assumptions on the middle surface.

For this reason, we solely focus on the residual error between the solution of Koiter's model and the averaged solution of the original three-dimensional model. We use conforming finite element to discretize the three components of the displacement(cf. , e.g., \cite{Shen2019}). Apart from the transverse component of the solution of Koiter's model, which is discretized by means of HCT triangles (cf., e.g., Chapter~6 of~\cite{PGCFEM}), all the other components of the solution of Koiter's model and of the solution of the three-dimensional model are approximated via a Lagrange finite element.

The domain $\omega$ is defined as follows 
\begin{align*} 
\omega:=\{(y_1,y_2)\in\mathbb{R}^2;(0,\pi)\times(0,1)\}, 
\end{align*} 
where $r=0.20 \meter,h=0.40 \meter$, and $\gamma_0$ is the region of the boundary at which the clamping occur. 

We use a portion of the cylindrical shell for numerical experiments and assume that the thickness of the shell is $2\varepsilon$. According to article \cite{Shen2019}, in curvilinear coordinates,the middle surface is given by the mapping $\bm{\theta}$ defined by 
\begin{align*} 
\bm{\theta}(y_1,y_2)=(r\cos y_1,r\sin y_1,hy_2)\quad\textup{ for all }(y_1,y_2)\in\overline{\omega}. 
\end{align*}

The second experiment we carry out concerns the displacement of a three-dimensional linearly elastic generalized membrane shell whose middle surface is a half sphere. We use the following values for the Lam\'e constants , Young's modulus, the Poission ratio and the applied body force density (cf., e.g., \cite{Shen2021}):
\begin{align*}
\nu&=0.45,\\
E&=5.4\times10^{6}\pascal,\\
\lambda&=1.68\times10^{7}\pascal,\\
\mu&=1.86\times10^{6}\pascal,\\
f^{i,\varepsilon}(x^\varepsilon)&=2.0\newton,\\
h^{i,\varepsilon}(x^\varepsilon)&=0\newton.
\end{align*}

In Table \ref{table:va2}, Err3DK indicates the residual error between the solution of Koiter's model and the solution of the 3D linearly elastic generalized membrane shell model. The consistent change in Table \ref{table:va2} means that as the thickness parameter $\varepsilon$ approaches zero, the residual error between the two aforementioned models tends to zero as well. 

\begin{figure}[h]
\centering
	\begin{subfigure}[b]{0.35\linewidth}
		\includegraphics[width=1.0\linewidth]{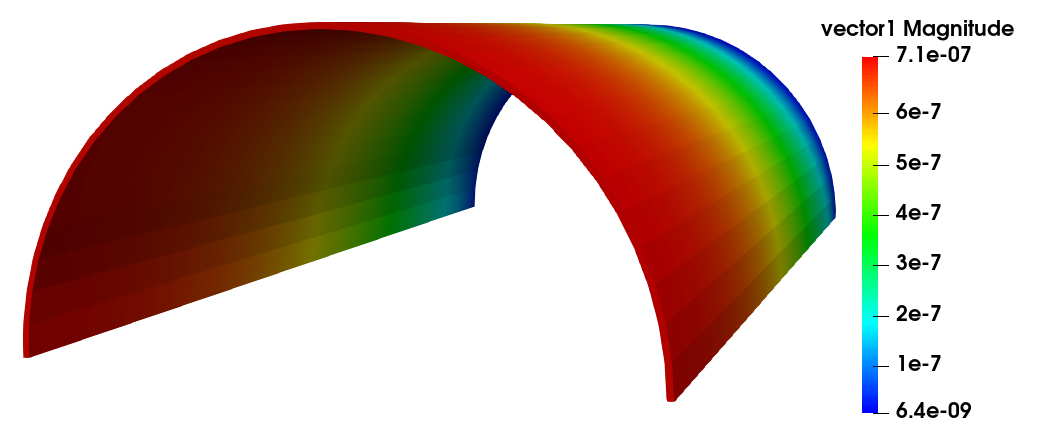}
		\subcaption{$\varepsilon=10^{-2}$}
	\end{subfigure}%
	\hspace{0.5cm}
	\begin{subfigure}[b]{0.35\linewidth}
		\includegraphics[width=1.0\linewidth]{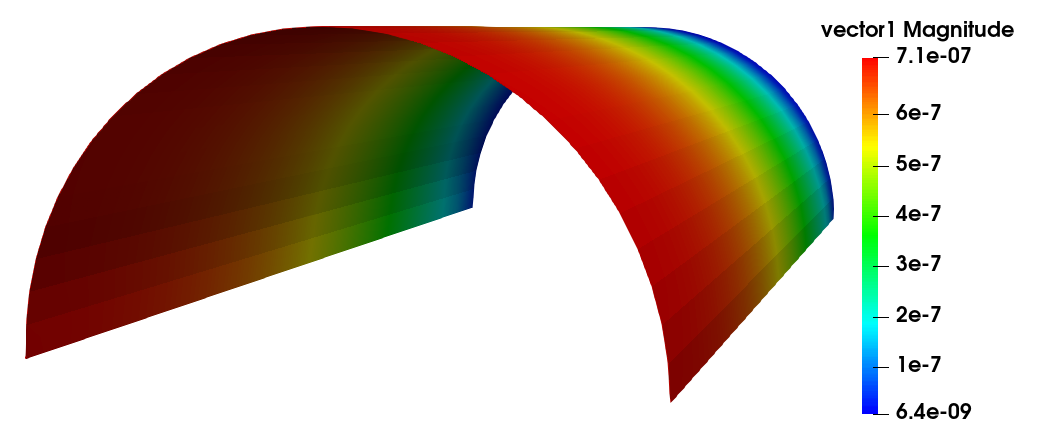}
		\subcaption{$\varepsilon=10^{-3}$}
	\end{subfigure}%
	\vspace{0.5cm}
	\begin{subfigure}[b]{0.35\linewidth}
		\includegraphics[width=1.0\linewidth]{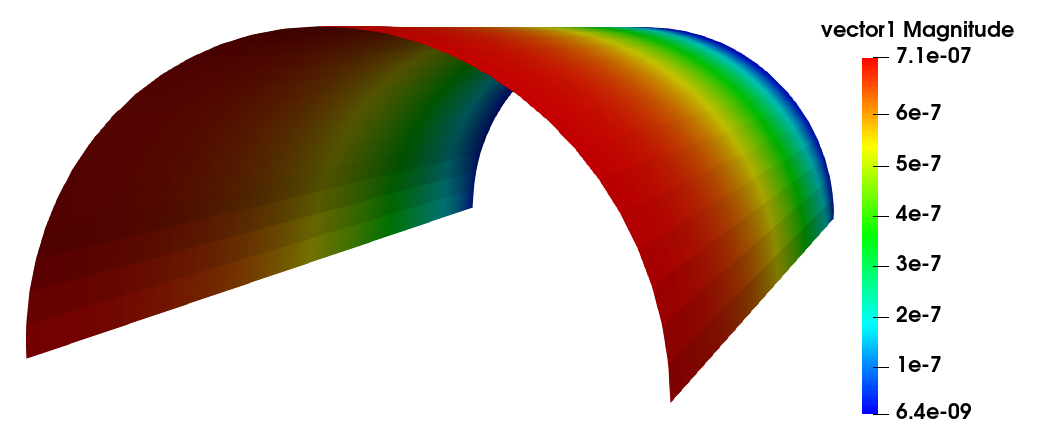}
		\subcaption{$\varepsilon=10^{-4}$}
	\end{subfigure}%
	\hspace{0.5cm}
	\begin{subfigure}[b]{0.35\linewidth}
		\includegraphics[width=1.0\linewidth]{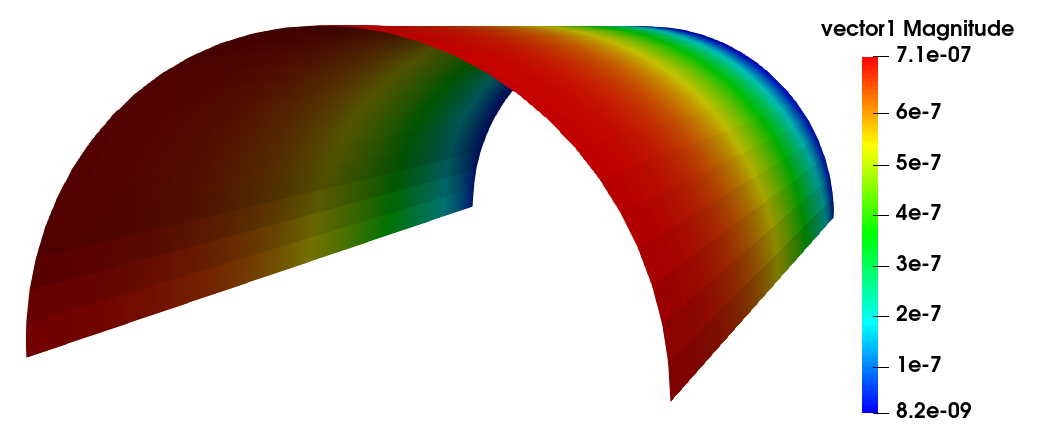}
		\subcaption{$\varepsilon=10^{-5}$}
	\end{subfigure}%
	\caption{Displacement of a three-dimensional linearly elastic generalized membrane shell whose middle surface is not elliptic. Different values of the parameter $\varepsilon$ are taken into account.}
	\label{CYGMS1}
\end{figure}

\begin{figure}[H]
\centering
	\begin{subfigure}[b]{0.35\linewidth}
		\includegraphics[width=1.0\linewidth]{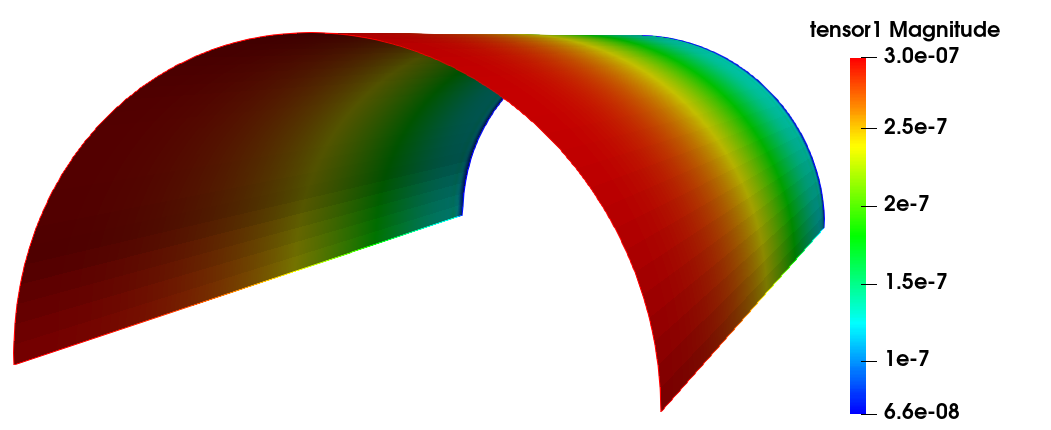}
		\subcaption{$\varepsilon=10^{-2}$}
	\end{subfigure}%
	\hspace{0.5cm}
	\begin{subfigure}[b]{0.35\linewidth}
		\includegraphics[width=1.0\linewidth]{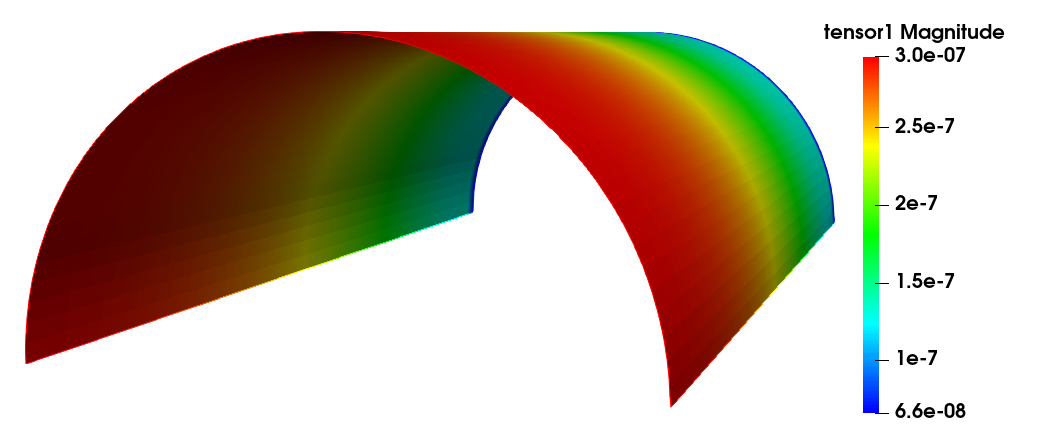}
		\subcaption{$\varepsilon=10^{-3}$}
	\end{subfigure}%
	\vspace{0.5cm}
	\begin{subfigure}[b]{0.35\linewidth}
		\includegraphics[width=1.0\linewidth]{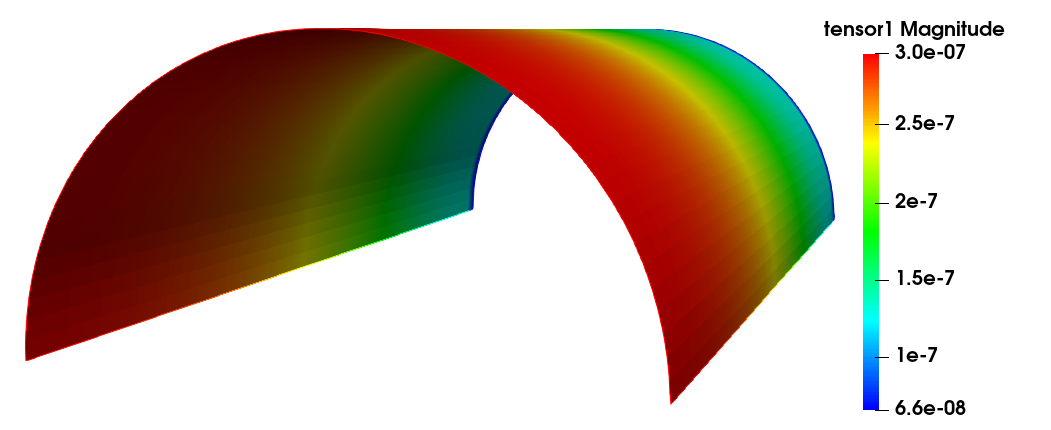}
		\subcaption{$\varepsilon=10^{-4}$}
	\end{subfigure}%
    \hspace{0.5cm}
	\begin{subfigure}[b]{0.35\linewidth}
		\includegraphics[width=1.0\linewidth]{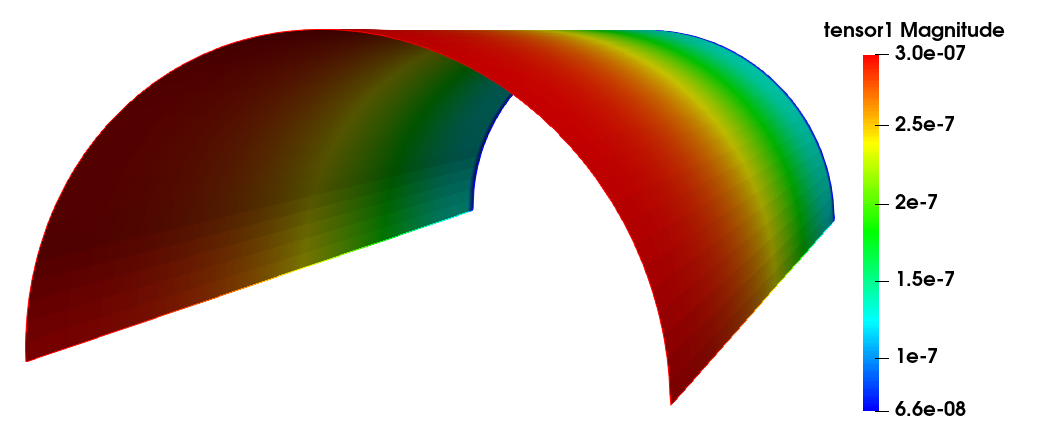}
		\subcaption{$\varepsilon=10^{-5}$}
	\end{subfigure}%
	\caption{Displacement of a two-dimensional Koiter's shell under the assumptions on the data characterizing a linearly elastic generalized membrane shell of the \emph{``first kind''}, viz., that the applied body forces are \emph{admissible}. Different values of the parameter $\varepsilon$ are taken onto account.}
	\label{CYGMS2}
\end{figure}

Figure \ref{CYGMS1} and \ref{CYGMS2} show the numerical results for the 3D linearly elastic generalized membrane shell model and Koiter's model. The entire bottom of each figure is in blue color, indicating that the entire boundary is clamped, whereas the biggest deformation happens on the parts in red color, i.e., the top of the cylindrical shell. It can be observed from the consistent shapes of the figures that the 3D linearly elastic elliptic membrane shell model and Koiter's model manifest the same asymptotic behavior as the thickness parameter $\varepsilon$ approaches zero.

\begin{table}[htbp]
\centering
\caption{Cylindrical generalized membrane shell result.}
\vspace{0.05cm}
\huge{
 \resizebox{\textwidth}{5mm}{\setlength{\tabcolsep}{3mm}{
  \begin{tabular}{c | c c c c c c c}
  \toprule[ 2pt]
    $\varepsilon$&1$\times10^{-2}$&1$\times10^{-3}$&1$\times10^{-4}$&1$\times10^{-5}$&1$\times10^{-6}$&1$\times10^{-7}$&1$\times10^{-8}$\\
    \hline
    Err3DK&3.10389e$-$07&1.00028e$-$07&3.16941e$-$08&1.00245e$-$08&3.1701e$-$09&1.00248e$-$09&3.17011e$-$10\\
\bottomrule[ 2pt]
  \end{tabular}\label{table:va2}
  }}}
\end{table}

\section{Numerical study of lineraly elastic flexural shells}\label{Sec:5}

In Section \ref{Sec:2}, the variational problem for ``general'' linearly elastic shells is considered by us. Based on the following definition, we will fix ourselves on a specific class of shells(originally proposed in \cite{Ciarlet19963}; see also \cite{Ciarlet2000}).

According to the assumptions set in Section 3, a \emph{linearly elastic flexural shell} will be studied in this Section, which is defined as follows: \emph{first}, $\emptyset \neq \gamma_0 \subset \gamma$, i.e., the homogeneous boundary condition of place is imposed over a \emph{non-zero area portion of the entire lateral face} $\gamma_0 \times \left[ - \varepsilon , \varepsilon \right]$ of the shell, and \emph{second}, the space
$$
\bm{V}_F(\omega):=\{\bm{\eta}=(\eta_i) \in H^1(\omega)\times H^1(\omega)\times H^2(\omega); \gamma_{\alpha \beta}(\bm{\eta})=0 \textup{ in }\omega
\textup{ and }\eta_i=\partial_{\nu}\eta_3=0 \textup{ on }\gamma_0\},
$$
contains non-zero functions, i.e., we have $\bm{V}_F(\omega)\neq\{\bm{0}\}$.

Let us define the sets
\[
\Omega := \omega \times \left] - 1, 1 \right[ \textup{ and } \Gamma_0:=\gamma_0 \times \left[ - 1, 1 \right] \textup{ and } \Gamma_{\pm}:=\omega \times \{\pm 1\},
\]
let $\partial_i := \partial/ \partial x_i$, and the generic point in the set $\overline{\Omega}$ is represented by $x = (x_i)$. With each point $x = (x_i) \in \overline{\Omega}$, we associate  the point $x^\varepsilon = (x^\varepsilon_i)$ defined by
\[
x^\varepsilon_\alpha := x_\alpha = y_\alpha \textup{ and } x^\varepsilon_3 := \varepsilon x_3,
\]
so that $\partial^\varepsilon_\alpha = \partial_\alpha$ and $\partial^\varepsilon_3 = \varepsilon^{-1} \partial_3$. 

For the unknown $\bm{u}^\varepsilon = (u^\varepsilon_i)$ and the vector fields $\bm{v}^\varepsilon = (v^\varepsilon_i)$ appearing in the formulation of the problem~\ref{problem0} corresponding to a linearly elastic flexural shell, we then associate the \emph{scaled unknown} $\bm{u} (\varepsilon) = (u_i(\varepsilon))$ and the \emph{scaled vector fields} $\bm{v} = (v_i)$ by letting
\[
u_i (\varepsilon) (x) := u^\varepsilon_i (x^\varepsilon) \text{ and } v_i(x) := v^\varepsilon_i (x^\varepsilon)
\]
at each $x\in \overline{\Omega}$. Finally, we \emph{assume} that there exist functions
, $f^i \in L^2(\Omega)$ and $h^i \in L^2(\Gamma_{+} \cup \Gamma_{-})$, are independent of $\varepsilon$,  which satisfy the following assumptions
\begin{equation}
	\label{dataFlex}
	\begin{aligned}
		f^{i, \varepsilon}(x^\varepsilon) &= \varepsilon^2 f^i(x) \textup{ for a.a. } x \in \Omega,\\
		h^{i, \varepsilon}(x^\varepsilon) &= \varepsilon^3 h^i(x) \textup{ for a.a. } x \in \Gamma_{+} \cup \Gamma_{-}.
	\end{aligned}
\end{equation}

Note that there is an implicit assumption for the Lam\'{e} constants that is independent of $\varepsilon$ assumed in the formulation of problem~\ref{problem0} in Section \ref{Sec:2}.

As a result of the scalings, the following identities hold (cf. section~\ref{Sec:3}):
\begin{align*}
	g(\varepsilon) (x) &:= g^\varepsilon (x^\varepsilon)
	\text{ and } A^{ijk\ell} (\varepsilon) (x) :=
	A^{ijk\ell, \varepsilon}(x^\varepsilon)
	\text{ at each } x\in \overline{\Omega}, \\
	e_{\alpha \| \beta} (\varepsilon; \bm{v}) &:= \frac12 (\partial_\beta v_\alpha + \partial_\alpha v_\beta ) - \Gamma^k_{\alpha \beta} (\varepsilon) v_k = e_{\beta \| \alpha} (\varepsilon; \bm{v}) , \\
	e_{\alpha \| 3} (\varepsilon; \bm{v}) = e_{3\| \alpha}(\varepsilon, \bm{v}) &:= \frac12 \left(\frac{1}{\varepsilon} \partial_3 v_\alpha + \partial_\alpha v_3 \right) - \Gamma^\sigma_{\alpha3} (\varepsilon) v_\sigma, \\
	e_{3\| 3} (\varepsilon ; \bm{v}) &:= \frac{1}{\varepsilon} \partial_3 v_3,
\end{align*}
where
\[
\Gamma^p_{ij} (\varepsilon) (x) := \Gamma^{p, \varepsilon}_{ij} (x^\varepsilon) \text{ at each } x\in \overline{\Omega}.
\]

We also let:
\[
\bm{g}_i (\varepsilon)(x) := g^\varepsilon_i (x^\varepsilon) \textup{ at each } x \in \Omega.
\]

\subsection{Asymptotic analysis}

Now we begin to state the convergence theorem for a family of linearly elastic flexural shells, that was originally proposed in~\cite{Ciarlet19963} (see also Theorem~6.2-1 of~\cite{Ciarlet2000}).

\begin{theorem}
\label{t:10}
Assume that $\bm{\theta} \in \mathcal{C}^2(\overline{\omega};\mathbb{R}^3)$. Consider a family of linearly elastic flexural shells with thickness $2\varepsilon$ approaching zero and with each having the same elliptic middle surface $\bm{\theta}(\overline{\omega})$, and let the assumptions on the data be as in Section~\ref{Sec:2} (for the Lam\'e constants) and~\eqref{dataFlex}.

For $\varepsilon>0$ sufficiently small, the solution of the associated scaled three-dimensional problems~\ref{problem0s} is represented by $\bm{u}(\varepsilon)$.

Let
\begin{align*}
a^{\alpha\beta\sigma\tau}&=\dfrac{4\lambda\mu}{\lambda+2\mu} a^{\alpha\beta} a^{\sigma\tau}+2\mu (a^{\alpha\sigma} a^{\beta\tau} + a^{\alpha\tau} a^{\beta\sigma}),\\
\gamma_{\alpha \beta}(\bm{\eta})&:=\frac12 (\partial_\beta \eta_\alpha + \partial_\alpha \eta_\beta ) - \Gamma^\sigma_{\alpha \beta} \eta_\sigma - b_{\alpha \beta} \eta_3,\\
\rho_{\alpha \beta}(\bm{\eta})&= \partial_{\alpha \beta}\eta_3 -\Gamma_{\alpha\beta}^\sigma \partial_\sigma\eta_3 -b_\alpha^\sigma b_{\sigma\beta}\eta_3+b_\alpha^\sigma(\partial_\beta \eta_\sigma-\Gamma_{\beta\sigma}^\tau\eta_\tau)\\
&\quad+b_\beta^\tau(\partial_\alpha\eta_\tau-\Gamma_{\alpha\tau}^\sigma\eta_\sigma)
+(\partial_\alpha b_\beta^\tau+\Gamma_{\alpha\sigma}^\tau b_\beta^\sigma-\Gamma_{\alpha\beta}^\sigma b_\sigma^\tau)\eta_\tau,\\
p^i&=\int_{-1}^{1} f^i \dd x_3 +h^i_{+}+h^i_{-}, \quad h^i_\pm=h^i(\cdot, \pm 1).
\end{align*}

Then there exist a functions on $\Gamma_0$, expressed as $u_\alpha \in H^1(\Omega)$, which satisfies $u_\alpha=0$, and a function $u_3 \in L^2(\Omega)$ such that:
\begin{align*}
	\bm{u}(\varepsilon) \to\bm{u}\quad&\textup{ in } \bm{H}^1(\Omega) \textup{ as }\varepsilon\to 0,\\
	\bm{u}=(u_i) &\textup{ is independent of the transverse variable } x_3.
\end{align*}

Furthermore, the \emph{average}
$$
\overline{\bm{u}}:=\dfrac{1}{2}\int_{-1}^{1} \bm{u} \dd x_3
$$
satisfies the following scaled variational two-dimensional problem $\mathcal{P}_F(\omega)$ of a linearly elastic flexural shell:
\begin{customprob}{$\mathcal{P}_F(\omega)$}\label{limitF}
	Find $\overline{\bm{u}} \in \bm{V}_F(\omega)=\{\bm{\eta}=(\eta_i) \in H^1(\omega)\times H^1(\omega)\times H^2(\omega);\\ \gamma_{\alpha \beta}(\bm{\eta})=0 \textup{ in }\omega
	\textup{ and }\eta_i=\partial_{\nu}\eta_3=0 \textup{ on }\gamma_0\}$ that satisfies the following variational equations:
	$$
	\int_{\omega} a^{\alpha\beta\sigma\tau} \rho_{\sigma\tau}(\overline{\bm{u}}) \rho_{\alpha \beta}(\bm{\eta}) \sqrt{a} \dd y = \int_{\omega} p^{i} \eta_i\sqrt{a} \dd y,
	$$
	for all $\bm{\eta} = (\eta_i) \in \bm{V}_F(\omega)$.
	\bqed	
\end{customprob}
The vector field $\bm{\zeta}$ is the one and only one which solves Problem~\ref{limitF}.
\qed
\end{theorem}

The results of Theorem \ref{t:10}, which can be applied to the solutions $\bm{u}(\varepsilon)$ of the scaled problem~\ref{problem0s}, still needs to be "de-scaled". It means that we need to change these results as ones about the \emph{unknown} $u^\varepsilon_i \bm{g}^{i, \varepsilon} : \overline{\Omega^\varepsilon} \to \mathbb{E}^3$, which represents the physical \emph{3D vector field} of the actual reference configuration of the shell. As we all know, this change is directly achieved by the introduction of the \emph{averages} $\displaystyle \frac{1}{2\varepsilon} \int^\varepsilon_{-\varepsilon} u^\varepsilon_i \bm{g}^{i, \varepsilon} \dd x^\varepsilon_3$ across the thickness of the shell. We will demonstrate the \emph{difference} (in terms of function spaces) between the asymptotic behaviours of the \emph{tangential} and \emph{normal} components of the displacement field of the middle surface of the shell.

\begin{theorem}[Theorem 6.4-1 in \cite{Ciarlet2000}] \label{t:11}
	Let the assumptions on the data be as in \emph{Section \ref{Sec:3}} and let the assumptions on the immersion $\bm{\theta} \in \mathcal{C}^3 (\overline{\omega}; \mathbb{E}^3)$ be as in \emph{Theorem \ref{t:10}}. Let $\bm{u}^\varepsilon = (u^\varepsilon_i) \in \bm{V}(\Omega^\varepsilon)$ denote for each $\varepsilon > 0$ the unique solution of the variational problem~\ref{problem0}, and let $\bm{\zeta} = (\zeta_i) \in \bm{V}_M(\omega)$ denote the unique solution for Problem~\ref{limitF}. Then
	\begin{align*}
		\frac{1}{2\varepsilon} \int^\varepsilon_{-\varepsilon} u^\varepsilon_\alpha \bm{g}^{\alpha , \varepsilon} \dd x^\varepsilon_3 & \to \zeta_\alpha \bm{a}^\alpha \textup{ in } \bm{H}^1 (\omega) \textup{ as } \varepsilon \to 0, \\
		\frac{1}{2\varepsilon} \int^\varepsilon_{-\varepsilon} u^\varepsilon_3 \bm{g}^{3 , \varepsilon} \dd x^\varepsilon_3 &\to \zeta_3 \bm{a}^3
		\textup{ in } \bm{H}^1 (\omega) \textup{ as } \varepsilon \to 0. 
	\end{align*}
	\qed
\end{theorem}

\subsection{Justification of Koiter's model for flexural shells}
\label{S:3:1}
The justification of Koiter's model for linearly elastic flexural shells was established by Ciarlet \& Lods in the paper~\cite[Theorem~2.2]{Ciarlet19962}.
The preamble is the same as in section~\ref{S:1:1}.

\begin{theorem}
	\label{t:k:3}
	Assume that $\bm{\theta} \in \mathcal{C}^2(\overline{\omega};\mathbb{R}^3)$. Consider a family of linearly elastic flexural shells with thickness $2\varepsilon$ approaching zero and with each having the same elliptic middle surface $\bm{\theta}(\overline{\omega})$, and let the assumptions on the data be as in Section~\ref{Sec:2} (for the Lam\'e constants) and~\eqref{dataFlex}. For each $\varepsilon>0$ let $\bm{u}^\varepsilon=(u_i^\varepsilon) \in \bm{H}^1(\Omega^\varepsilon)$ and $\bm{\zeta}_K^\varepsilon=(\zeta_{i,K}^\varepsilon)$ denote the solutions to the three-dimensional problem~\ref{problem0} and the two-dimensional problem~\ref{Koiter}.
	Let also
	$$
	\bm{\zeta} =(\zeta_i) \in \bm{V}_M(\omega)= H^1(\omega) \times H^1(\omega) \times H^2(\omega),
	$$
	denote the solution to the two-dimensional scaled variational problem~\ref{limitM}, solution which is independent of $\varepsilon$. Then, the following convergences hold:
	\begin{align*}
		\dfrac{1}{2\varepsilon} \int_{-\varepsilon}^{\varepsilon} u^\varepsilon_\alpha \bm{g}^{\alpha,\varepsilon} \dd x_3^\varepsilon \to \zeta_\alpha \bm{a}^\alpha &\textup{ in } \bm{H}^1(\omega) \textup{ as } \varepsilon \to 0,\\
		\zeta_{\alpha,K}^\varepsilon \bm{a}^\alpha \to \zeta_\alpha \bm{a}^\alpha &\textup{ in } \bm{H}^1(\omega) \textup{ as } \varepsilon \to 0,
	\end{align*}
	and
	\begin{align*}
		\dfrac{1}{2\varepsilon} \int_{-\varepsilon}^{\varepsilon} u^\varepsilon_3 \bm{g}^{3,\varepsilon} \dd x_3^\varepsilon \to \zeta_3 \bm{a}^3 &\textup{ in } \bm{H}^1(\omega) \textup{ as } \varepsilon \to 0,\\
		\zeta_{3,K}^\varepsilon \bm{a}^3 \to \zeta_3 \bm{a}^3 &\textup{ in } \bm{H}^2(\omega) \textup{ as } \varepsilon \to 0.
	\end{align*}
	\qed
\end{theorem}

Note that the first and the third convergence in Theorem~\ref{t:k:3} have been recovered in Theorem~\ref{t:11}.

\subsection{Numerical results}
We verify the convergences announced in Theorem~\ref{t:11} through a series of \emph{ad hoc} numerical experiments carried out by implementing the finite element method (cf., e.g., \cite{PGCFEM}) via the software FreeFem~\cite{Hecht2012}. The results are visualised in ParaView~\cite{Ahrens2005}. Specifically, We use conforming finite element to discretize the three components of the displacement(cf., e.g., \cite{ShenXue2021}). Apart from the transverse component of the solution of Koiter's model, which is discretized by means of HCT triangles (cf., e.g., Chapter~6 of~\cite{PGCFEM}), all the other components of the solution of Koiter's model, the solution of the 3D model and the solution of the 2D limit model are approximated via a Lagrange finite element.

The domain $\omega$ is defined as follows 
\begin{align*} 
\omega:=\{(y_1,y_2)\in\mathbb{R}^2;(0,\pi)\times (0.4,1)\}, 
\end{align*} 
where $b=0.20 \meter,c=0.40 \meter$. 

We use a portion of the conical shell for numerical experiments and assume that the thickness of the shell is $2\varepsilon$. According to article \cite{ShenXue2021}, in curvilinear coordinates, the middle surface is given by the mapping $\bm{\theta}$ defined by
\begin{align*} 
\bm{\theta}(y_1,y_2)=(by_2 \cos y_1,by_2 \sin y_1,cy_2)\quad\textup{ for all } (y_1,y_2)\in\overline{\omega}. 
\end{align*} 

The third experiment we carry out concerns the displacement of a three-dimensional linearly elastic flexural shell whose middle surface is a portion of a conical shell. We use the following values for the Lam\'e constants, Young's modulus, the Poisson ratio and the applied body force density (cf., e.g., \cite{Shen2021}):
\begin{align*}
\nu&=0.45,\\
E&=5.40\times10^{6}\pascal,\\
\lambda&=1.68\times10^{7}\pascal,\\
\mu&=1.86\times10^{6}\pascal,\\
f^{i,\varepsilon}(x^\varepsilon)&=8.00\times10^{1}\newton,\\
h^{i,\varepsilon}(x^\varepsilon)&=0\newton.
\end{align*}

In Table \ref{table:va3}, ErrLK indicates the residual error between the solution for the 2D limit model and Koiter's model; Err3DL indicates the residual error between the solution of the 3D linearly elastic flexural shell model and the 2D limit model; Err3DK is the error of the solution between Koiter's model and the 3D linearly elastic flexural shell model. The consistent change in Table \ref{table:va3} means that as the thickness parameter $\varepsilon$ approaches zero, the residual errors between the different models tend to zero as well.

\begin{figure}[!http]
	\centering
	\includegraphics[width=0.85\linewidth]{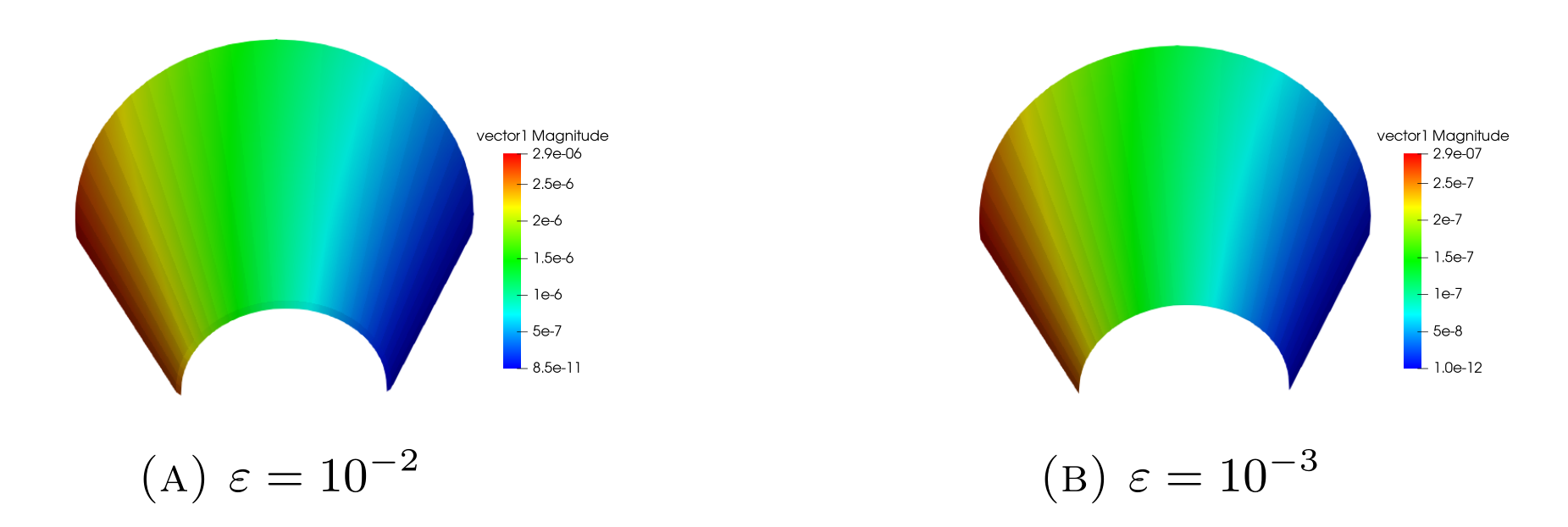}
	\includegraphics[width=0.85\linewidth]{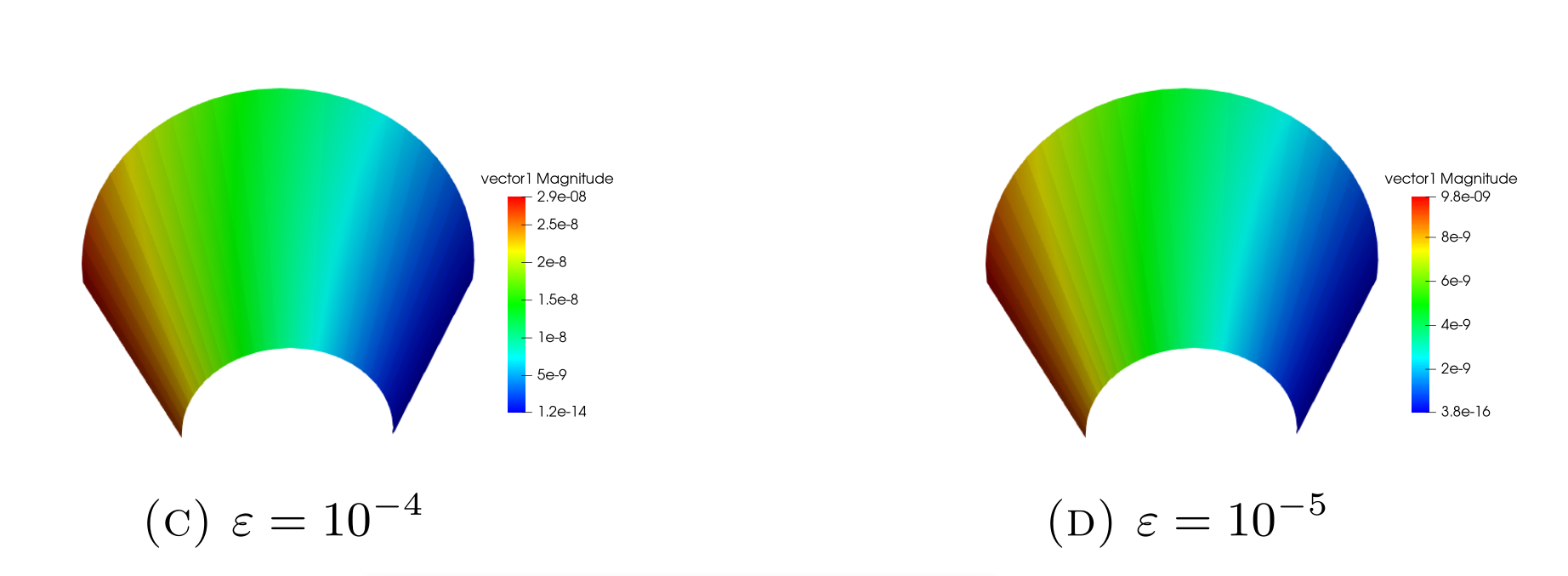}
	\caption{Displacement of a three-dimensional linearly elastic flexural shell whose middle surface is a portion of a cone. Different values of the parameter $\varepsilon$ are taken into account.}
	\label{COFS1}
\end{figure}

\begin{figure}[!http]
	\centering
	\includegraphics[width=0.85\linewidth]{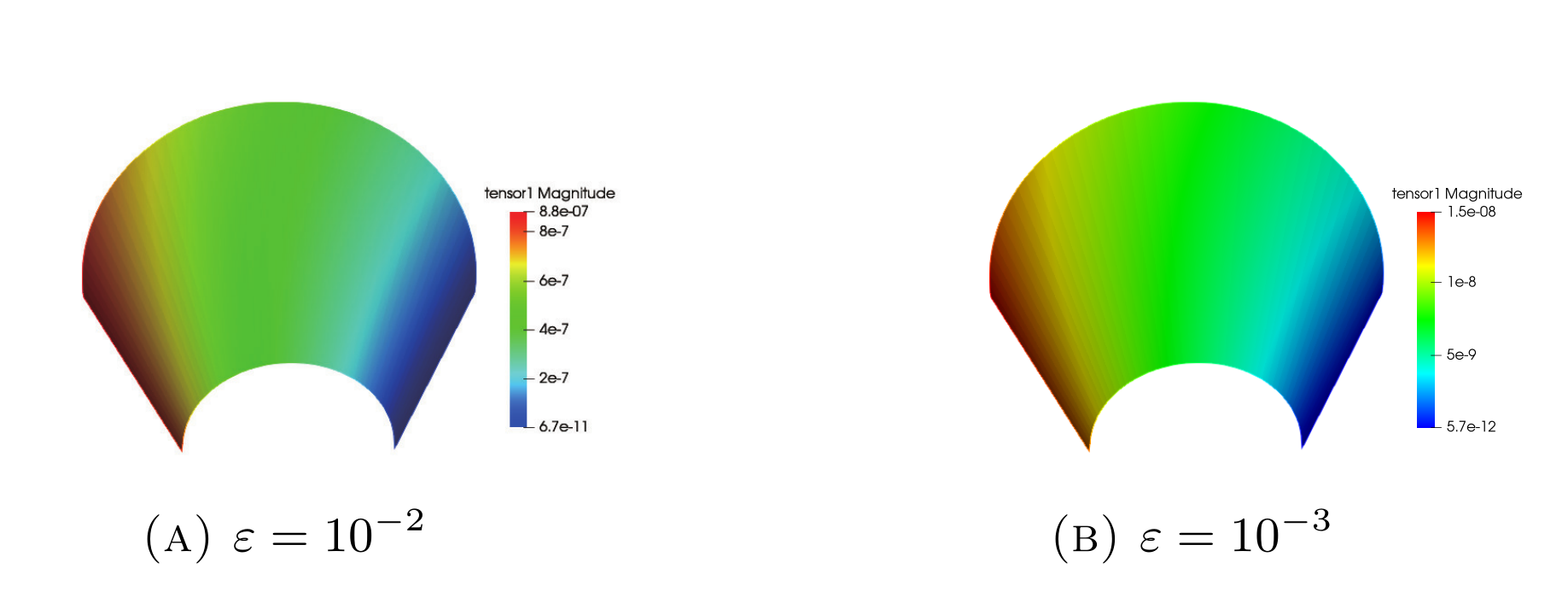}
	\includegraphics[width=0.85\linewidth]{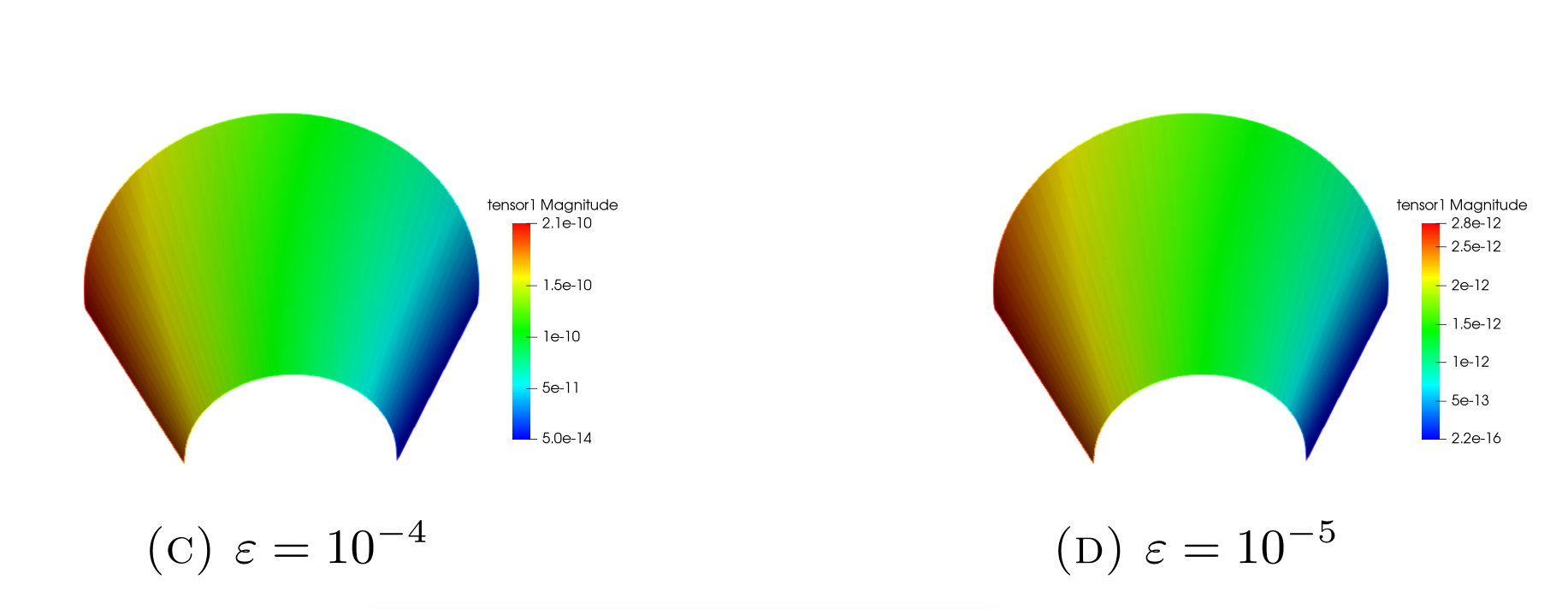}
	\caption{Displacement of a two-dimensional Koiter's shell under the assumption on the data~\eqref{dataFlex} characterizing linearly elastic flexural shells. Different values of the parameter $\varepsilon$ are taken into account.}
	\label{COFS2}
\end{figure}

Figure \ref{COFS1} and \ref{COFS2} show the numerical results of the 3D linearly elastic flexural shell model and Koiter's model. The right generatrix of each figure is in blue, indicating that the entire boundary is clamped, whereas the biggest deformation happens on the parts in red color, i.e., the left generatrix of the conical shell. It can be observed from the consistent shapes of the figures that the 3D linearly elastic flexural shell model and Koiter's model manifest the same asymptotic behavior as the thickness parameter $\varepsilon$ approaches zero.

\begin{table}[htbp]
\centering
\caption{Conical flexural shell result.}
\vspace{0.05cm}
\huge{
 \resizebox{\textwidth}{16mm}{\setlength{\tabcolsep}{20mm}{
  \begin{tabular}{c c c c}
  \toprule[2pt]
  \specialrule{0em}{3pt}{3pt}
   $\varepsilon$ &ErrLK  & Err3DK & Err3DL\\
   \midrule[1pt]
   \specialrule{0em}{3pt}{3pt}
    5.02655e$-$01 & 1.16982e$-$06 & 1.14079e$-$06 & 2.90302e$-$08  \\
    5.02655e$-$02 & 2.01333e$-$08 & 1.98439e$-$08  & 2.89427e$-$10   \\
    5.02655e$-$03 & 2.97578e$-$10 & 2.94712e$-$10  & 2.86643e$-$12   \\
    5.02655e$-$04 & 4.26798e$-$12 & 4.16869e$-$12  & 9.92919e$-$14   \\
    5.02655e$-$05 & 4.71892e$-$14 & 4.71861e$-$14 & 3.05265e$-$18 \\
    5.02655e$-$06 & 4.72919e$-$16 & 4.72918e$-$16 & 5.39086e-$-$22 \\
\bottomrule[2pt]
  \end{tabular}\label{table:va3}
}}}

\end{table}

\newpage

\section{Conclusions and future works}
\label{Sec:6}

In this paper we verified the justification of Koiter's model in the static case via \emph{ad hoc} numerical experiments. We showed that, for all the three main types of linearly elastic shells (viz., linearly elastic elliptic membrane shells, linearly elastic generalized membrane shells of the first kind, and linearly elastic flexural shells), the average across the thickness of the solution for the three-dimensional model asymptotically behaves like the solution of Koiter's model as the thickness parameter $\varepsilon$ approaches zero.

In future works, we will verify the justification of Koiter's model via numerical examples in the time-dependent case, that was analytically addressed in the papers~\cite{Xiao2001,Xiao2001flex}, and in the case where the displacement of the shell under consideration is also affected by a thermal source, viz., the shell is thermoelastic (cf., e.g., \cite{CaoRodCasRos2021} and~\cite{Pie2021}).

\section*{Acknowledgements}

The second author are greatly indebted to Professor Philippe G. Ciarlet and Professor Roger M. T\'{e}mam for their encouragement and guidance. 

The third author is greatly indebted to Professor Philippe G. Ciarlet for his encouragement and guidance. 

This paper is supported by the Ky and Yu-Fen Fan Fund Travel Grant from the AMS (IU Award number 44-294-36 with the Simons Foundation) and the National Natural Science Foundation of China (NSFC.11971379).	

\bibliographystyle{abbrvnat} 
\bibliography{references.bib}	

\end{document}